# DIRECTIONAL GAUSSIAN HYPERGEOMETRIC BETA DISTRIBUTIONS AND THEIR USES IN CONTAMINATED BINARY SAMPLING

BEN O'NEILL[*] *ACIL Allen*[**]

CANBERRA, AUSTRALIA

ORCID: 0000-0001-6899-0483

## Abstract

We examine the Gaussian hypergeometric beta distribution and look at the effect of having an additional term in the density kernel relative to the standard beta distribution. We reparameterise and classify this distribution into left and right directional variants using parameters that give a simple and symmetrical representation of the directional push/pull from this additional term in the density kernel. We examine the properties of the directional variants and their uses in contaminated binary sampling using Bayesian inference. We find that the Gaussian hypergeometric beta distribution arises as the appropriate posterior distribution for inference in certain kinds of contaminated binary models and that the directional parameterisation aids in representation of the resulting Bayesian models. We derive a broad range of properties and computational methods for the directional variants of the distribution.

[*] E-mail address: ben.oneill@hotmail.com.
[**] Level 6, 54 Marcus Clarke St, Canberra ACT 2601, Australia.

## 1. The Gaussian hypergeometric beta (HyperBeta) distribution

The beta distribution is widely used in statistics to model unknown probability parameters and it arises commonly in Bayesian models for binary data. In particular, since the beta distribution is the conjugate prior for the probability parameter in the binomial distribution it arises as both the prior and posterior form in the commonly used beta-binomial model. The beta-binomial model is well-known, but useful variations of this model involving generalisations of the beta distribution are not widely known.

Extensions of the beta distribution in statistical literature have been used to capture conjugate priors for other forms of Bayesian model. These have sometimes involved adding additional monotonic terms to the kernel of the beta density which "push" or "pull" the bulk of the density to the left or right of the support in various ways. A useful generalisation of this kind is the Gaussian hypergeometric beta (HyperBeta) distribution, introduced in Armero and Bayarri (1994) and Gordy (1998). In this paper we examine this distributional family in the context of a simple model for binary sampling subject to a controlled contamination mechanism. We will show that the HyperBeta distribution arises as the conjugate prior for the likelihood function arising from contaminated binary sampling, and this forms a HyperBeta-binomial model that generalises the beta-binomial model for the non-contaminated case.

Whilst the HyperBeta distribution is mathematically suited to this task, its parameterisation in standard form is somewhat awkward in regard to representing its "push" and "pull" properties relative to the left and right boundaries of the support. We find that it is useful to reparameterise the HyperBeta distribution and split it into two directional variants that give a simple and symmetrical representation of its push/pull properties relative to the left or right of the support. This leads to an alternative parameterised form that is mathematically equivalent to the usual form, but simpler to interpret and apply to the contaminated binomial model (and possibly other problems). We will introduce this alternative parameterisation/classification of the HyperBeta distribution and examine various properties of the distribution, including its application to the contaminated binomial model. Our examination of the HyperBeta distribution covers core properties, including probability functions, shape, moments and computation methods. We will first show the density function for our parameterisation of the distribution and explain its derivation in the contaminated binary model. We then examine the shape of the distribution, including the monotonicity and quasi-concavity properties of the directional variants, and how these depend on the parameters. We establish the moments of the distribution and some other moment-related properties and we derive the loglikelihood and score functions. Many of these results involve integrals that do not have a closed form, so we also examine computational methods to compute various functions relating to this distribution.

Before proposing our preferred parameterisation, we first introduce the standard form for the Gaussian hypergeometric beta (HyperBeta) distribution. An examination of the properties of this distribution can be found in Nadarajah and Kebe (2024), and later we will conduct some similar analysis working with our own parameterisation. The standard form for the density function for the HyperBeta distribution is fully characterised by its density kernel:[1]

$$\text{HyperBeta}(x|\alpha,\beta,\gamma,\delta) \propto x^{\alpha-1}(1-x)^{\beta-1}(1+\delta x)^{\gamma} \qquad 0 \leq x \leq 1,$$

where $\alpha > 0$ and $\beta > 0$ are shape parameters analogous to those in the beta distribution and $-\infty < \gamma < \infty$ and $\delta > -1$ are additional parameters that generalise the beta distribution. The density function is obtained by adding the relevant scaling constant, giving the form:

$$\begin{aligned}\text{HyperBeta}(x|\alpha,\beta,\gamma,\delta) &= \frac{x^{\alpha-1}(1-x)^{\beta-1}(1+\delta x)^{\gamma}}{\int_0^1 x^{\alpha-1}(1-x)^{\beta-1}(1+\delta x)^{\gamma} dx} \\ &= \frac{x^{\alpha-1}(1-x)^{\beta-1}(1+\delta x)^{\gamma}}{\mathrm{B}(\alpha,\beta)\ {}_2F_1(-\gamma,\alpha,\alpha+\beta|-\delta)}.\end{aligned}$$

The scaling constant for the density function is an integral that can be represented using the beta function and the Gaussian hypergeometric function, given respectively by:[2]

$$\mathrm{B}(\alpha,\beta) = \int_0^1 x^{\alpha-1}(1-x)^{\beta-1} dx\,,$$

$${}_2F_1(-\gamma,\alpha,\alpha+\beta|-\delta) = \sum_{n=0}^{\infty} \frac{(\gamma)_n(\alpha+n-1)_n}{(\alpha+\beta+n-1)_n} \cdot \frac{\delta^n}{n!}\,,$$

where $(r)_n \equiv \prod_{i=0}^{n-1}(r-i)$ are the falling factorials. The distribution takes its name from the presence of the Gaussian hypergeometric function in the adjustment to the beta distribution. The Gaussian hypergeometric function provides the solution to a second-order linear ordinary differential equation (Temme 1996, pp. 107-132; Andrews, Askey and Roy 1999, pp. 61-123).

---

[1] Nadarajah and Kebe (2024) use a parameterisation using $c = -\gamma$ and $d = \delta$. We have made a minor variation here, partly to simplify later analysis and partly to avoid using the notation $d$ for a parameter (since resulting multiples of the form $dx$ can easily be confused with differentials in later analysis).

[2] The Gaussian hypergeometric function is sometimes presented in terms of rising factorials; we have used falling factorials here to simplify terms that include negative coefficients.

To understand our reason for reparameterising the distribution, it is useful to first examine the directional properties of the additional term in the kernel of the HyperBeta distribution relative to the standard beta distribution. Relative to the beta distribution, we have the additional term $(1+\delta x)^\gamma$ in the density kernel of the HyperBeta distribution. This is a monotonic function of $x$ that moves the bulk of the density towards one of the boundaries of the support, depending on the parameters $\gamma$ and $\delta$. There are three directional cases of interest shown in Table 1.

**TABLE 1: Directional properties of the Gaussian hypergeometric distribution**

| Distribution Type | Condition | Directional shift (relative to beta distribution) |
|---|---|---|
| Beta | $\delta = 0$ | No shift in distribution (relative to beta distribution) |
| Left-HyperBeta | $\delta < 0$ | Directional shift towards the **left** of the support as the parameter $\gamma$ increases. For all values $-\infty < \gamma < \acute{\gamma} < \infty$ we have the monotone likelihood ratio: $$\frac{d}{dx}\frac{\text{HyperBeta}(x\vert\alpha,\beta,\acute{\gamma},\delta)}{\text{HyperBeta}(x\vert\alpha,\beta,\gamma,\delta)} < 0 \qquad 0 < x < 1.$$ The additional term in the kernel pulls the density away from the left of the support if $\gamma < 0$ or pushes the density towards the left of the support if $\gamma > 0$. |
| Right-HyperBeta | $\delta > 0$ | Directional shift towards the **right** of the support as the parameter $\gamma$ increases. For all values $-\infty < \gamma < \acute{\gamma} < \infty$ we have the monotone likelihood ratio: $$\frac{d}{dx}\frac{\text{HyperBeta}(x\vert\alpha,\beta,\acute{\gamma},\delta)}{\text{HyperBeta}(x\vert\alpha,\beta,\gamma,\delta)} > 0 \qquad 0 < x < 1.$$ The additional term in the kernel pulls the density away from the right of the support if $\gamma < 0$ or pushes the density towards the right of the support if $\gamma > 0$. |

In the standard parameterisation of the HyperBeta distribution, the density obeys the reflection property $\text{HyperBeta}(x|\alpha,\beta,\gamma,\delta) = \text{HyperBeta}(1-x|\beta,\alpha,\gamma,\acute{\delta})$ with the adjusted parameter value $\acute{\delta} = -\delta/(1+\delta)$, which means that the HyperBeta family captures reflection involving equivalent forms of pushing/pulling to both the left and right of the support, albeit in a form that uses a non-symmetric parameterisation. Rather than using the standard parameterisation, we find it is useful to parameterise the distribution to use a *symmetric representation* of these directional shifts relative to the left and right of the support. This is done by separating the distribution into two classes based on the sign of the parameter $\delta$, yielding left-HyperBeta and right-HyperBeta distributions, with an accompanying transformation of the parameter $\delta$.

## 2. The directional HyperBeta distributions

Now that we have examined the directional properties of the HyperBeta distribution relative to the beta distribution, we can stipulate a useful reparameterised form that allows symmetric characterisation of the directional variants. Our parameter transformation is designed to give left and right versions of the distribution that capture the relevant push/pull in a simple and symmetrical way. To do this, we define a binary indicator $\kappa = \mathbb{I}(\delta > 0)$ to represent the left and right cases of the distribution and we define a new push proportion parameter $0 \leq \phi \leq 1$ via the correspondence:

$$\delta = \begin{cases} -\phi & \text{if } \kappa = 0, \\ \phi/(1-\phi) & \text{if } \kappa = 1. \end{cases}$$

This parameterisation gives the following density kernels in the left and right cases:

$$\begin{array}{lll} \kappa = 0 & (1+\delta x)^{\gamma} \propto (1 - x\phi)^{\gamma} & \text{LHyperBeta} \\ \kappa = 1 & (1+\delta x)^{\gamma} \propto (1 - \phi + x\phi)^{\gamma} & \text{RHyperBeta} \end{array}$$

Our reparameterisation separates the HyperBeta distribution into left and right directional variants, which will be shown to obey a simpler version of the reflection principle. We define the reparameterised form of the distribution here by stipulating its probability density.

**DEFINITION 1 (Gaussian hypergeometric beta distribution):** This is a continuous univariate distribution with support on the unit interval. It has density over $0 \leq x \leq 1$ given by:

$$\text{HyperBeta}(x|\alpha,\beta,\gamma,\phi,\kappa) = \frac{x^{\alpha-1}(1-x)^{\beta-1}(1-x\phi)^{(1-\kappa)\gamma}(1-\phi+x\phi)^{\kappa\gamma}}{\int_0^1 x^{\alpha-1}(1-x)^{\beta-1}(1-x\phi)^{(1-\kappa)\gamma}(1-\phi+x\phi)^{\kappa\gamma}dx}.$$

where $\alpha, \beta > 0$ are the **shape parameters**, $-\infty < \gamma < \infty$ is the **push-intensity parameter**, $0 \leq \phi \leq 1$ is the **push-proportion parameter** and $\kappa \in \{0, 1\}$ is the **push-direction** ($\kappa = 0$ for the left-push and $\kappa = 1$ for the right-push). The density can also be written as:[3]

$$\text{HyperBeta}(x|\alpha,\beta,\gamma,\phi,\kappa) = \frac{x^{\alpha-1}(1-x)^{\beta-1}(1-x\phi)^{(1-\kappa)\gamma}(1-\phi+x\phi)^{\kappa\gamma}}{\text{B}(\alpha,\beta)\ {}_2F_1(-\gamma,\alpha^{1-\kappa}\beta^{\kappa},\alpha+\beta|\phi)}.$$

This distribution is a generalisation of the standard beta distribution with an additional term in the density kernel reflecting a push/pull to the left or right side of the support. □

[3] Here B refers to the beta function and ${}_2F_1$ refers to the hypergeometric function. The equivalence of the two forms for the scaling constant in the density function is an expression of Euler's integral formula for the hypergeometric function (see e.g., Bailey 1935, pp. 4-5).

Our reparameterisation of the HyperBeta distribution includes several special cases that reduce down to the standard beta distribution. In particular, if the push-intensity parameter is zero, we obtain the special case:

$$\gamma = 0 \qquad \text{HyperBeta}(x|\alpha, \beta, 0, \phi, \kappa) = \text{Beta}(x|\alpha, \beta)\,,$$

which reflects there being no push in the distribution. Moreover, at the extremes of the range of the push-proportion parameter we have:

$$\phi = 0 \qquad \text{HyperBeta}(x|\alpha, \beta, \gamma, 0, \kappa) = \text{Beta}(x|\alpha, \beta)\,,$$

$$\phi = 1 \qquad \text{HyperBeta}(x|\alpha, \beta, \gamma, 1, \kappa) = \text{Beta}(x|\alpha + \kappa\gamma, \beta + (1-\kappa)\gamma)\,,$$

which reflects that these extreme pushes can be absorbed into the standard beta form. Our reparameterisation preserves the ability to represent any Gaussian hypergeometric distribution.

The above definition shows the general form of the HyperBeta distribution, but it is simplest to view this family of distributions as the union of two subfamilies which are the left-HyperBeta distribution ($\kappa = 0$) and right-HyperBeta distribution ($\kappa = 1$). These directional variants of the HyperBeta distribution have respective density functions given by:

$$\text{LHyperBeta}(x|\alpha, \beta, \gamma, \phi) = \text{HyperBeta}(x|\alpha, \beta, \gamma, \phi, 0) = \frac{x^{\alpha-1}(1-x)^{\beta-1}(1-x\phi)^{\gamma}}{\text{B}(\alpha, \beta)\ {}_2F_1(-\gamma, \alpha, \alpha+\beta|\phi)},$$

$$\text{RHyperBeta}(x|\alpha, \beta, \gamma, \phi) = \text{HyperBeta}(x|\alpha, \beta, \gamma, \phi, 1) = \frac{x^{\alpha-1}(1-x)^{\beta-1}(1-\phi+x\phi)^{\gamma}}{\text{B}(\alpha, \beta)\ {}_2F_1(-\gamma, \beta, \alpha+\beta|\phi)}.$$

From these forms we can see that the density kernels for the beta and HyperBeta distributions have a similar form, with the latter being a slight generalisation of the former:

$$\text{Beta}(x|\alpha, \beta) \propto x^{\alpha-1}(1-x)^{\beta-1},$$

$$\text{LHyperBeta}(x|\alpha, \beta, \gamma, \phi) \propto x^{\alpha-1}(1-x)^{\beta-1}(1-x\phi)^{\gamma},$$

$$\text{RHyperBeta}(x|\alpha, \beta, \gamma, \phi) \propto x^{\alpha-1}(1-x)^{\beta-1}(1-\phi+x\phi)^{\gamma}.$$

The form of the HyperBeta density is similar to the standard beta density, but with an additional push/pull term. For the left-HyperBeta distribution we have the additional term $(1-x\phi)^{\gamma}$ which is downward sloping for $\phi > 0$ and $\gamma > 0$ and upward sloping for $\phi > 0$ and $\gamma < 0$. For the right-HyperBeta distribution we have the additional term $(1-\phi+x\phi)^{\gamma}$ which is upward sloping for $\phi > 0$ and $\gamma > 0$ and downward sloping for $\phi > 0$ and $\gamma < 0$. From the reparameterisation we now obtain the simplified **reflection equations**:

$$\text{LHyperBeta}(x|\alpha, \beta, \gamma, \phi) = \text{RHyperBeta}(1-x|\beta, \alpha, \gamma, \phi)\,,$$

$$\text{RHyperBeta}(x|\alpha, \beta, \gamma, \phi) = \text{LHyperBeta}(1-x|\beta, \alpha, \gamma, \phi)\,.$$

In Theorems 1-2 below we establish that greater values of $\phi$ or $\gamma$ push the density further to the left or right boundaries of the support (for the left-HyperBeta and right-HyperBeta cases respectively) in terms of first-order stochastic dominance.

**THEOREM 1 (Stochastic dominance with respect to push proportion):** Suppose we take the push-proportions $0 \leq \phi_0 < \phi_1 \leq 1$ and define the random variables:

$$X_0 \sim \text{HyperBeta}(\alpha, \beta, \gamma, \phi_0, \kappa), \qquad X_1 \sim \text{HyperBeta}(\alpha, \beta, \gamma, \phi_1, \kappa) .$$

Then we have the following stochastic dominance results depending on the parameters:

- If $\kappa = 0$ and $\gamma > 0$ then $X_0$ stochastically dominates $X_1$ (written as $X_0 \succ X_1$);
- If $\kappa = 0$ and $\gamma < 0$ then $X_1$ stochastically dominates $X_0$ (written as $X_1 \succ X_0$);
- If $\kappa = 1$ and $\gamma > 0$ then $X_1$ stochastically dominates $X_0$ (written as $X_1 \succ X_0$);
- If $\kappa = 1$ and $\gamma < 0$ then $X_0$ stochastically dominates $X_1$ (written as $X_0 \succ X_1$).

**THEOREM 2 (Stochastic dominance with respect to push intensity):** Suppose we take the push intensities $-\infty < \gamma_0 < \gamma_1 < \infty$ with $0 < \phi < 1$ and define the random variables:

$$X_0 \sim \text{HyperBeta}(x|\alpha, \beta, \gamma_0, \phi, \kappa) \qquad X_1 \sim \text{HyperBeta}(x|\alpha, \beta, \gamma_1, \phi, \kappa) .$$

Then we have the following stochastic dominance results depending on the parameters:

- If $\kappa = 0$ then $X_0$ stochastically dominates $X_1$ (written as $X_0 \succ X_1$);
- If $\kappa = 1$ then $X_1$ stochastically dominates $X_0$ (written as $X_1 \succ X_0$).

Before proceeding to a more detailed examination of the directional variants of the HyperBeta distribution we first show a simple problem involving contaminated binary sampling where this distribution arises. Our contaminated binary sampling problem involves a situation where a researcher is interested in determining the prevalence of a "sensitive" (binary) characteristic pertaining to people in some population of interest, where the characteristic is of such a sensitive nature that subjects may lie about the characteristic top the researcher, either due to social desirability bias and/or other negative consequence from honest disclosure. Confidentiality protections may go some way towards alleviating this risk, but for highly sensitive characteristics it may be insufficient to yield honest disclosure, due to the subject's perception of a risk of a confidentiality breach. In such cases, an available protocol to reduce subject concerns is for the researcher to impose "controlled contamination" on the observation of the sensitive characteristic, to obscure any answers that would disclose the presence of the sensitive characteristic at issue (even in the event that there is a confidentiality breach in the responses). To do this, the researcher has the subject add a controlled random event to each response (e.g., a die roll) and solicits the information in such a way that disclosure of the sensitive characteristic is "contaminated" by the contamination event (which is only observed by the subject and not the researcher). This controlled contamination event serves to make it impossible for the researcher to determine whether an individual subject has the sensitive characteristic of interest, which also means that the subject does not rely on the confidentiality protections in the research. The known probability of contamination still allows valid statistical inference of the overall presence of the sensitive characteristic in the population.

This idea can be better understood by showing the contaminated binary sampling model and then giving further explanation of the experimental protocols at issue. Suppose we observe a sequence of binary outcomes formed as indicators of the conjunction of two sets of conditions (a "primary condition" for the characteristic of interest and a "contaminating condition"). Since the contaminating condition is under the control of the researcher, we will assume that the probability of the contaminating condition is known and fixed over the trials, and we are seeking to make an inference about the unknown probability of the primary condition. We encapsulate this situation by taking two independent sequences of Bernoulli random variables (we will call these the primary and contaminating sequences respectively) given by:

$$\tilde{X}_1, \tilde{X}_2, \tilde{X}_3, \dots \sim \text{IID Bern}(\theta) \qquad Y_1, Y_2, Y_3, \dots \sim \text{IID Bern}(\phi)\,.$$

The random variable $\tilde{X}_i$ is the indicator for the primary condition with unknown probability parameter $\theta$ which we want to infer from data, and the random variable $Y_i$ is the indicator for the contaminating condition which has (known) probability parameter $\phi$. In this problem we observe the product values $X_i = \tilde{X}_i Y_i$ from the sequence:

$$X_1, X_2, X_3, \dots \sim \text{IID Bern}(\theta\phi)\,.$$

In this model, the contaminating events $Y_1, Y_2, Y_3, \dots$ are experimental protocols imposed by the researcher and so the probability $\phi$ is a known control parameter in the experiment. It can be varied by the researcher to obscure a sensitive characteristic for an individual subject. In the above case, the characteristic $\tilde{X}_i = 0$ is the "sensitive" result (e.g., some characteristic that is socially undesirable) but the corresponding observation $X_i = 0$ does not imply the presence of that characteristic. If the researcher sets $\phi$ to be sufficiently small then an observation $X_i = 0$ still gives the subject "plausible deniability" that $\tilde{X}_i = 0$ even if the result is disclosed (e.g., due to a confidentiality breach in the research).

**REMARK:** Because of the contamination of the primary sequence (yielding a product value), one might be tempted to view this observation mechanism as being similar to incomplete binary data with missing values, by considering $Y_i = 0$ to denote a kind of "missingness" at random. However, in the contaminated binary model these "missing" values are indistinguishable from true zeroes for the primary sequence —i.e., we cannot distinguish the cases $\tilde{X}_i = 0, Y_i = 1$, $\tilde{X}_i = 1, Y_i = 0$ or $\tilde{X}_i = 0, Y_i = 0$, which would be distinguishable in standard missing data problems (e.g., with the missing data outcome denoted as a separate outcome $X_i = \text{NA}$). Unlike in a binary missing data problem, here there is contamination occurring where we can only see whether or not the binary values are *both unity*. The assumption that we know the probability parameter $\phi$ for the contaminating sequence is crucial to the model because the product $\theta\phi$ is the minimal sufficient parameter (O'Neill 2005) — it is not possible to identify either $\theta$ or $\phi$ if both are unknown. □

The left-HyperBeta distribution arises from this model using a simple Bayesian approach with a conjugate prior for the unknown probability parameter $\theta$. To see this, suppose we observe sample values $\boldsymbol{x}_n = (x_1, x_2, \dots, x_n)$ and remind ourselves that the parameter $0 \leq \phi \leq 1$ is taken to be known. Taking $\dot{x}_n = \sum_{i=1}^n x_i$ to be the sample sum we get the sampling density:

$$f(\boldsymbol{x}_n|\theta) = \prod_{i=1}^{n} (\theta\phi)^{x_i}(1-\theta\phi)^{1-x_i} = (\theta\phi)^{\dot{x}_n}(1-\theta\phi)^{n-\dot{x}_n},$$

which establishes that $\dot{x}_n$ is a sufficient statistic for $\theta$. (In fact, it can easily be shown that it is minimal sufficient.) Using the prior $\theta \sim \text{LHyperBeta}(\alpha, \beta, \gamma, \phi)$ for the unknown probability parameter (with a beta prior as a special case) we get the posterior density:

$$\begin{aligned}
\pi(\theta|\boldsymbol{x}_n) &\propto f(\boldsymbol{x}_n|\theta) \cdot \pi(\theta) \\
&\propto (\theta\phi)^{\dot{x}_n}(1-\theta\phi)^{n-\dot{x}_n} \cdot \theta^{\alpha-1}(1-\theta)^{\beta-1}(1-\theta\phi)^{\gamma} \\
&\propto \theta^{\alpha+\dot{x}_n-1}(1-\theta)^{\beta-1}(1-\theta\phi)^{\gamma+n-\dot{x}_n} \\
&\propto \text{LHyperbeta}(\theta|\alpha+\dot{x}_n, \beta, \gamma+n-\dot{x}_n, \phi)\,.
\end{aligned}$$

Note that the left-HyperBeta posterior arises even if we use a standard beta prior (with scale parameters $\alpha$ and $\beta$) since this is equivalent to using a left-HyperBeta prior with zero push-intensity. This means that the left-HyperBeta distribution arises as a conjugate prior under the standard treatment of an unknown probability parameter when we observe binary outcomes subject to the contamination in the model. This contaminated binary model is summarised in Table 2, showing the occurrence of the left-HyperBeta distribution as a conjugate distributional family.

| **TABLE 2: Contaminated Binary Model** | |
|---|---|
| We observe data from the sequence: $$X_1, X_2, X_3, \ldots \sim \text{IID Bern}(\theta\phi),$$ where $0 \leq \phi \leq 1$ is a known value and $0 \leq \theta \leq 1$ is an unknown parameter of interest. (Note that $\phi$ must be known or the parameter $\theta$ is unidentifiable in the model.) | |
| **Sample vector** | $\boldsymbol{x}_n = (x_1, x_2, \ldots, x_n)$ |
| **Minimal sufficient statistic** | $\dot{x}_n = \sum_{i=1}^n x_i$ |
| **Sampling density** | $f(\boldsymbol{x}_n \vert \theta) = (\theta\phi)^{\dot{x}_n}(1-\theta\phi)^{n-\dot{x}_n}$ |
| **Likelihood function** | $L_{\boldsymbol{x}_n}(\theta) \propto \theta^{\dot{x}_n}(1-\theta\phi)^{n-\dot{x}_n}$ |
| **Conjugate prior and posterior** | $\pi(\theta) = \text{LHyperBeta}(\theta \vert \alpha, \beta, \gamma, \phi)$<br>$\pi(\theta \vert \boldsymbol{x}_n) = \text{LHyperBeta}(\theta \vert \alpha + \dot{x}_n, \beta, \gamma + n - \dot{x}_n, \phi)$ |
| **Special cases** | In the case where $\gamma = 0$ we get the prior:<br>$\pi(\theta) = \text{LHyperBeta}(\theta \vert \alpha, \beta, 0, \phi) = \text{Beta}(\theta \vert \alpha, \beta)$<br>The beta distribution is commonly used to model an unknown probability parameter so this special case is of interest. |

The above model gives rise to the left-HyperBeta distribution as a conjugate distribution. but a simple variation on this model gives rise to the right-HyperBeta distribution. This occurs when contamination operates on the absence of the primary condition instead of its presence. In this case we instead observe the product values $X_i = (1 - \tilde{X}_i)Y_i$ from the sequence:

$$X_1, X_2, X_3, \ldots \sim \text{IID Bern}((1-\theta)\phi).$$

The right-HyperBeta distribution arises from this model in an analogous way to the previous case. This alternative contaminated binary model gives us the sampling density:

$$f(\boldsymbol{x}_n \vert \theta) = \prod_{i=1}^n ((1-\theta)\phi)^{x_i}(1-\phi+\theta\phi)^{1-x_i} = ((1-\theta)\phi)^{\dot{x}_n}(1-\phi+\theta\phi)^{n-\dot{x}_n}.$$

Using the prior $\theta \sim \text{RHyperBeta}(\alpha, \beta, \gamma, \phi)$ for the unknown probability parameter (with a beta prior as a special case) we get the posterior density:

$$\begin{aligned}\pi(\theta \vert \boldsymbol{x}_n) &\propto f(\boldsymbol{x}_n \vert \theta) \cdot \pi(\theta) \\ &\propto ((1-\theta)\phi)^{\dot{x}_n}(1-\phi+\theta\phi)^{n-\dot{x}_n} \cdot \theta^{\alpha-1}(1-\theta)^{\beta-1}(1-\phi+\theta\phi)^{\gamma} \\ &\propto \theta^{\alpha-1}(1-\theta)^{\beta+\dot{x}_n-1}(1-\phi+\theta\phi)^{\gamma+n-\dot{x}_n} \\ &\propto \text{RHyperBeta}(\theta \vert \alpha, \beta + \dot{x}_n, \gamma + n - \dot{x}_n, \phi).\end{aligned}$$

This means that the right-HyperBeta distribution arises as a conjugate prior under the standard treatment of an unknown probability parameter when we observe binary outcomes subject to the contamination in the model. This model is summarised in Table 3 below, showing the occurrence of the right-HyperBeta distribution as a conjugate distributional family.

| **TABLE 3: Contaminated Binary Model (Alternate)** | |
|---|---|
| We observe data from the sequence: $X_1, X_2, X_3, \dots \sim \text{IID Bern}((1-\theta)\phi)$, where $0 \leq \phi \leq 1$ is a known value and $0 \leq \theta \leq 1$ is an unknown parameter of interest. (Note that $\phi$ must be known or the parameter $\theta$ is unidentifiable in the model.) | |
| **Sample vector** | $\boldsymbol{x}_n = (x_1, x_2, \dots, x_n)$ |
| **Minimal sufficient statistic** | $\dot{x}_n = \sum_{i=1}^n x_i$ |
| **Sampling density** | $f(\boldsymbol{x}_n \mid \theta) = ((1-\theta)\phi)^{\dot{x}_n}(1-\phi+\theta\phi)^{n-\dot{x}_n}$ |
| **Likelihood function** | $L_{\boldsymbol{x}_n}(\theta) \propto (1-\theta)^{\dot{x}_n}(1-\phi+\theta\phi)^{n-\dot{x}_n}$ |
| **Conjugate prior and posterior** | $\pi(\theta) = \text{RHyperBeta}(\theta \mid \alpha, \beta, \gamma, \phi)$<br>$\pi(\theta \mid \boldsymbol{x}_n) = \text{RHyperBeta}(\theta \mid \alpha, \beta + \dot{x}_n, \gamma + n - \dot{x}_n, \phi)$ |
| **Special cases** | In the case where $\gamma = 0$ we get the prior:<br>$\pi(\theta) = \text{RHyperBeta}(\theta \mid \alpha, \beta, 0, \phi) = \text{Beta}(\theta \mid \alpha, \beta)$<br>The beta distribution is commonly used to model an unknown probability parameter so this special case is of interest. |

**REMARK:** The HyperBeta distribution is actually more general than needed to get a conjugate distribution for the contaminated binary models. In the model in Table 2, the term $(1-\theta)^{\beta-1}$ does not need to be included in the density kernel to maintain conjugacy with the likelihood function in this model. Similarly, in the model in Table 3, the term $\theta^{\alpha-1}$ does not need to be included in the density kernel to maintain conjugacy with the likelihood function in this model. Notwithstanding this extra generality, it is still useful to use the HyperBeta distribution since it is a conjugate prior in both model forms and it is also a generalisation of the beta distribution, which the minimal conjugate forms are not. These broader conjugate models extend the standard beta-binomial model that occurs in Bayesian analysis of non-contaminated binary sampling. □

We will examine practical applications of the contaminated binary model in a later section of this paper. For now it is worth noting that this is an elementary model which involves a simple variation of standard observation of IID binary random variables. This simple model gives rise to the directional HyperBeta distributions as conjugate forms in Bayesian analysis.

**3. Shape of the directional HyperBeta distributions**

We have now seen how the directional variants of the HyperBeta distribution arise in the contaminated binary model and we have also seen that —relative to the standard beta density— these involve pushing the density in a particular direction in a manner determined by the push-intensity $\gamma$ and push-proportion $\phi$. The next property we examine in this section is the shape of the density function. To get a good idea of the shape we - examine the parameter conditions leading to various monotonicity and quasi-concavity/quasi-convexity results for the density functions. We will also compare these conditions to the case of a standard beta distribution to get an idea of the effect of the "push" in the directional variants of the HyperBeta distribution.

To facilitate this analysis, let $\ell_x(\alpha, \beta, \gamma, \phi, \kappa) \equiv \log \text{HyperBeta}(x|\alpha, \beta, \gamma, \phi, \kappa)$ denote the log-density, which has derivatives given by:

$$\frac{d}{dx}\ell_x(\alpha, \beta, \gamma, \phi, \kappa) = \frac{\alpha - 1}{x} - \frac{\beta - 1}{1 - x} - \frac{(1 - \kappa)\gamma\phi}{1 - x\phi} + \frac{\kappa\gamma\phi}{1 - \phi + x\phi},$$

$$\frac{d^2}{dx^2}\ell_x(\alpha, \beta, \gamma, \phi, \kappa) = -\frac{\alpha - 1}{x^2} - \frac{\beta - 1}{(1 - x)^2} - \frac{(1 - \kappa)\gamma\phi^2}{(1 - x\phi)^2} - \frac{\kappa\gamma\phi^2}{(1 - \phi + x\phi)^2}.$$

The first derivative of the log-density tells us about monotonicity of the density function and the second derivative tells us about quasi-concavity/convexity of the density function. In both cases the third/fourth term in the expression is the contribution of the "push" part of the density, one of which will be non-zero if $\gamma \neq 0$ and $\phi > 0$. Assuming that these conditions hold, the push contribution in the first derivative is more negative for higher values of the push intensity in the left-HyperBeta distribution and more positive for higher values of the push intensity in the right-HyperBeta distribution. This means that higher push intensity gives more left-push in the left-HyperBeta distribution and more right-push in the right-HyperBeta distribution. The push contribution in the second derivative is negative for positive push intensity and more positive for negative push intensity, which means that this makes the density more concave (less convex) or more convex (less concave) than it would otherwise be in these two cases.

To obtain formal results for the shape of the density function it is useful to look more deeply at the direction of the slope by looking at the sign of the derivative of the log-density. Letting $\acute{\alpha} = \alpha - 1$ and $\acute{\beta} = \beta - 1$ we can write the derivative of the log-density in expanded form as:

$$\frac{d}{dx}\ell_x(\alpha,\beta,\gamma,\phi,\kappa) = \frac{Q_\kappa(x|\alpha,\beta,\gamma,\phi)}{x(1-x)(1-x\phi)^{1-\kappa}(1-\phi+x\phi)^\kappa},$$

where the numerator is the (left or right) quadratic function:

$$Q_0(x|\alpha,\beta,\gamma,\phi) \equiv \acute{\alpha} - (\acute{\alpha} + \acute{\beta} + \acute{\alpha}\phi + \gamma\phi)x + \phi(\acute{\alpha} + \acute{\beta} + \gamma)x^2,$$

$$Q_1(x|\alpha,\beta,\gamma,\phi) \equiv \acute{\alpha}(1-\phi) - [\acute{\alpha}(1-2\phi) + \acute{\beta}(1-\phi) - \gamma\phi]x - \phi(\acute{\alpha} + \acute{\beta} + \gamma)x^2.$$

The term in the denominator of the derivative is strictly positive over the interior of the support so the sign of the density function is determined by the sign of the quadratic function in the numerator (either $Q_0$ or $Q_1$). Denoting the roots of the first quadratic by $r_{0,0}, r_{0,1} \in \mathbb{C}$ and the roots of the second quadratic by $r_{1,0}, r_{1,1} \in \mathbb{C}$ gives the condensed forms for these quadratics:

$$Q_0(x|\alpha,\beta,\gamma,\phi) = \phi(\acute{\alpha} + \acute{\beta} + \gamma)(x - r_{0,0})(x - r_{0,1}),$$

$$Q_1(x|\alpha,\beta,\gamma,\phi) = -\phi(\acute{\alpha} + \acute{\beta} + \gamma)(x - r_{1,0})(x - r_{1,1}).$$

We also have the following values at the boundaries of the support:

$$Q_0(0|\alpha,\beta,\gamma,\phi) = \acute{\alpha} \qquad Q_1(0|\alpha,\beta,\gamma,\phi) = \acute{\alpha}(1-\phi),$$

$$Q_0(1|\alpha,\beta,\gamma,\phi) = -\acute{\beta}(1-\phi) \qquad Q_1(1|\alpha,\beta,\gamma,\phi) = -\acute{\beta}.$$

We will denote the direction of the slope of the HyperBeta density by:

$$\text{DIR}_\kappa(x|\alpha,\beta,\gamma,\phi) \equiv \text{sgn}\frac{d}{dx}\text{HyperBeta}(x|\alpha,\beta,\gamma,\phi,\kappa) = \text{sgn}\, Q_\kappa(x|\alpha,\beta,\gamma,\phi)\,.$$

The shape of the density function is determined by the signs at the boundaries of the support and the placement of the roots — in particular, whether the roots are real values and whether they fall within the support. If there are no real roots in the interior of the support then the density function has a single sign over the interior of the support and so it is monotonic (with direction that can be determined by the sign of the quadratic at the boundaries of the support). If there is one real root in the interior of the support then the density function changes sign only once over the interior of the support so it is either quasi-concave or quasi-convex (which can be determined by the sign of the quadratic at the boundaries of the support). Finally, if there are two distinct real roots in the interior of the support then the density function changes sign twice over the interior of the support, which means that it is neither quasi-concave nor quasi-convex. In Theorems 3-4 below we formalise this result, which leads to various monotonicity and quasi-concavity/quasi-convexity results depending on the parameters.

**THEOREM 3 (Slope of density – reflection equation):** We have the reflection equation:

$$\text{DIR}_\kappa(x|\alpha,\beta,\gamma,\phi) = -\text{DIR}_{1-\kappa}(1-x|\beta,\alpha,\gamma,\phi).$$

**THEOREM 4 (Slope of density – left and right HyperBeta):** The direction of the slope of the density in the left and right cases is given respectively by:

$$\text{DIR}_0(x|\alpha,\beta,\gamma,\phi) = \text{sgn}(\alpha+\beta+\gamma-2) \times \text{sgn}(x - r_{0,0}) \times \text{sgn}(x - r_{0,1}) \,,$$

$$\text{DIR}_1(x|\alpha,\beta,\gamma,\phi) = -\,\text{sgn}(\alpha+\beta+\gamma-2) \times \text{sgn}(1 - x - r_{0,0}) \times \text{sgn}(1 - x - r_{0,1}) \,.$$

Theorem 3 gives us a simple principle for the reflection of the direction of the slope of the density in the left and right cases for the HyperBeta density. Theorem 4 then establishes the direction of the slope of the density over the support. Similar results for the shape of the Gaussian hypergeometric distribution (under its standard parameterisation) can be found in Nadarajah and Kebe (2024) (pp. 3-6). This gives a broad range of specific monotonicity and/or quasi-concavity/convexity results depending on the parameters. It is notable that there are cases where both roots fall within the unit interval so that the density is neither monotonic nor quasi-concave or quasi-convex. In general there may be up to two critical points of the density within the support and so it is possible to have shapes involving one-and-a-half "waves" that can result in various combinations of local modes and anti-modes.

It is useful to get a sense of the strength of the "push" that occurs in the directional HyperBeta density through an appropriate visualisation of how the density changes as we change the relevant parameters affecting this push. In Figures 1A-1B below we show the left-Hyperbeta density and look at the strength of the "left push" as we vary either the push intensity parameter or the push proportion parameter. In Figure 1A we use fixed parameters $\alpha$, $\beta$ and $\phi$ but with an increasing value of the push-intensity parameter $\gamma$. This figure shows the effect of adding the "left push" to a beta distribution with increasing intensity. We see that as the push-intensity increases we see the density function is pushed more to the left of the support. In this case we have used parameters that yield a quasi-concave density. In Figure 1B we use fixed parameters $\alpha$, $\beta$ and $\gamma$ but with an increasing value of the push-proportion parameter $\phi$. This figure shows the effect of adding the "left push" to a beta distribution with increasing proportion. We see that as the push-intensity increases we see the density function is pushed more to the left of the support. We have again used parameters that yield a quasi-concave density.

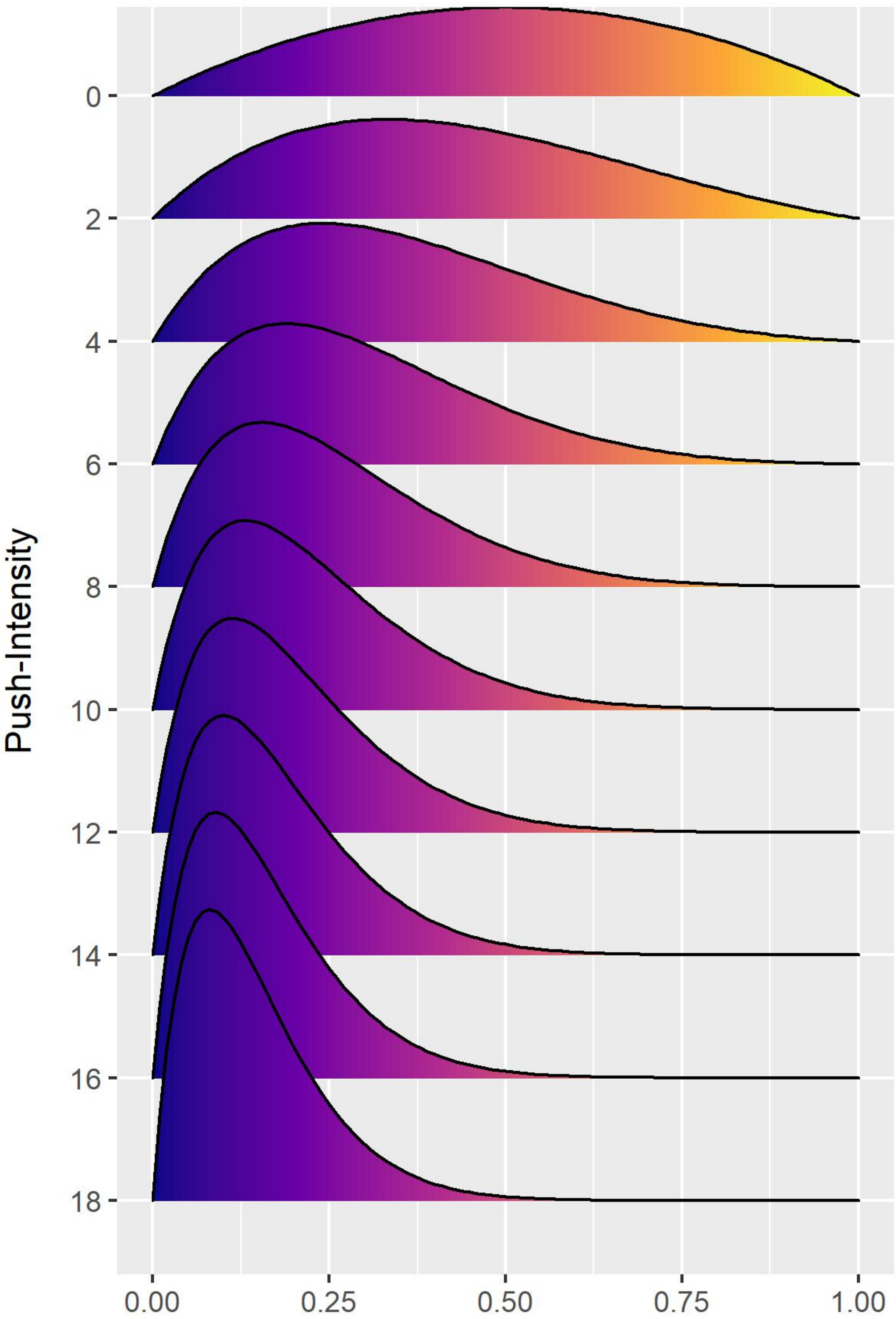

**FIGURE 1A:** Density function for the left-Hyperbeta distribution shown for various positive values of the push-intensity parameter (Other parameters fixed at $\alpha = 2$, $\beta = 2$ and $\phi = 0.6$)

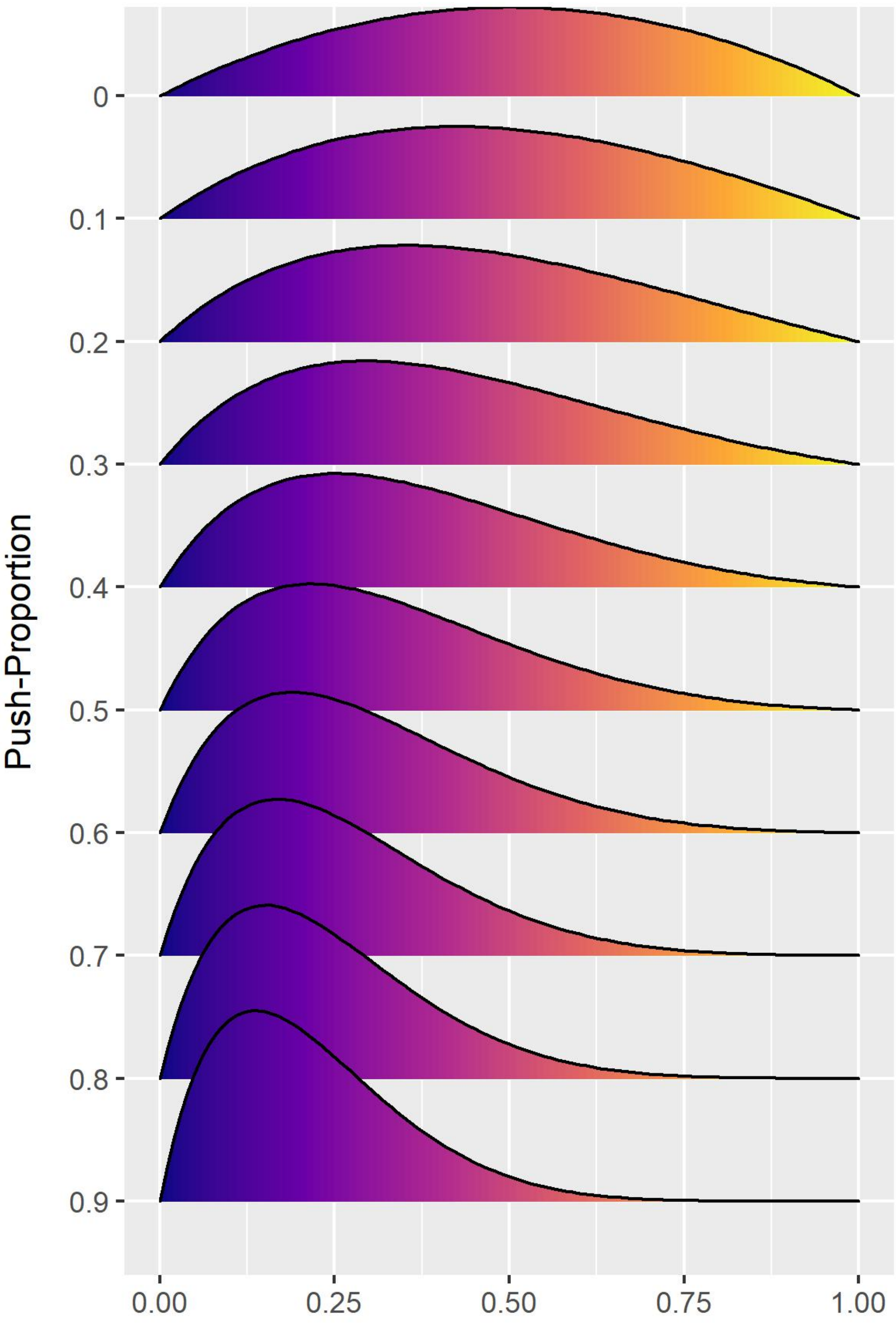

**FIGURE 1B:** Density function for the left-Hyperbeta distribution shown for various positive values of the push-proportion parameter (Other parameters fixed at $\alpha = 2$, $\beta = 2$ and $\gamma = 6$)

## 4. Moments and asymptotic properties of the directional HyperBeta distributions

The above results give a summary of the shape of the HyperBeta distribution in all its various cases. To understand the location, scale, etc., for the distribution we can examine its moments. In Theorem 6 below we show the raw moments of the HyperBeta distribution. As with the density function these do not have a closed form and they involve integration of the density kernel, which is equivalent to computation of the beta and hypergeometric functions. The mean and variance formulae are also shown; higher-order central moments can be written in terms of the raw moments shown, but their expressions are not particularly illuminating so they are omitted.

**REMARK (Moment reflections):** The moments of the two directions of the distribution are related through the reflection $\text{LHyperBeta}(x|\alpha,\beta,\gamma,\phi) = \text{RHyperBeta}(1-x|\beta,\alpha,\gamma,\phi)$ so that if we take $X_L \sim \text{LHyperBeta}(x|\alpha,\beta,\gamma,\phi)$ and $X_R \sim \text{LHyperBeta}(x|\beta,\alpha,\gamma,\phi)$ then we have the moment relations $\mathbb{E}(X_L^k) = \mathbb{E}((1-X_R)^k)$ and $\mathbb{E}(X_R^k) = \mathbb{E}((1-X_L)^k)$. □

**THEOREM 5A (Raw moments of left-HyperBeta beta distribution):** Suppose we have:

$$X \sim \text{LHyperBeta}(x|\alpha,\beta,\gamma,\phi).$$

This random variable has raw moments:

$$\mathbb{E}(X^k) = \frac{\text{B}(\alpha+k,\beta)}{\text{B}(\alpha,\beta)} \cdot \frac{{}_2F_1(-\gamma,\alpha+k,\alpha+\beta+k;\phi)}{{}_2F_1(-\gamma,\alpha,\alpha+\beta;\phi)} \qquad \text{for all } k = 1,2,3,\dots.$$

These can be written in alternative form as:

$$\mathbb{E}(X^k) = \frac{\int_0^1 x^{\alpha+k-1}(1-x)^{\beta-1}(1-x\phi)^\gamma dx}{\int_0^1 x^{\alpha-1}(1-x)^{\beta-1}(1-x\phi)^\gamma dx}.$$

**COROLLARY (Mean and variance):** If $X \sim \text{LHyperBeta}(x|\alpha,\beta,\gamma,\phi)$ then we have:

$$\mathbb{E}(X) = \frac{\text{B}(\alpha+1,\beta)}{\text{B}(\alpha,\beta)} \cdot \frac{{}_2F_1(-\gamma,\alpha+1,\alpha+\beta+1;\phi)}{{}_2F_1(-\gamma,\alpha,\alpha+\beta;\phi)},$$

$$\mathbb{V}(X) = \frac{\text{B}(\alpha+2,\beta)}{\text{B}(\alpha,\beta)} \cdot \frac{{}_2F_1(-\gamma,\alpha+2,\alpha+\beta+2;\phi)}{{}_2F_1(-\gamma,\alpha,\alpha+\beta;\phi)} - \left(\frac{\text{B}(\alpha+1,\beta)}{\text{B}(\alpha,\beta)} \cdot \frac{{}_2F_1(-\gamma,\alpha+1,\alpha+\beta+1;\phi)}{{}_2F_1(-\gamma,\alpha,\alpha+\beta;\phi)}\right)^2.$$

**THEOREM 5B (Raw moments of right-HyperBeta distribution):** Suppose we have:

$$X \sim \text{RHyperBeta}(x|\alpha,\beta,\gamma,\phi).$$

This random variable has raw moments:

$$\mathbb{E}(X^k) = \frac{\text{B}(\alpha+k,\beta)}{\text{B}(\alpha,\beta)} \cdot \frac{{}_2F_1(-\gamma,\beta,\alpha+\beta+k;\phi)}{{}_2F_1(-\gamma,\beta,\alpha+\beta;\phi)} \qquad \text{for all } k = 1,2,3,\dots.$$

These can be written in alternative form as:

$$\mathbb{E}(X^k) = \frac{\int_0^1 x^{\alpha+k-1}(1-x)^{\beta-1}(1-\phi+x\phi)^\gamma dx}{\int_0^1 x^{\alpha-1}(1-x)^{\beta-1}(1-\phi+x\phi)^\gamma dx}.$$

**COROLLARY (Mean and variance):** If $X \sim \text{RHyperBeta}(x|\alpha,\beta,\gamma,\phi)$ then we have:

$$\mathbb{E}(X) = \frac{\text{B}(\alpha+1,\beta)}{\text{B}(\alpha,\beta)} \cdot \frac{{}_2F_1(-\gamma,\beta,\alpha+\beta+1;\phi)}{{}_2F_1(-\gamma,\beta,\alpha+\beta;\phi)},$$

$$\mathbb{V}(X) = \frac{\text{B}(\alpha+2,\beta)}{\text{B}(\alpha,\beta)} \cdot \frac{{}_2F_1(-\gamma,\beta,\alpha+\beta+2;\phi)}{{}_2F_1(-\gamma,\beta,\alpha+\beta;\phi)} - \left(\frac{\text{B}(\alpha+1,\beta)}{\text{B}(\alpha,\beta)} \cdot \frac{{}_2F_1(-\gamma,\beta,\alpha+\beta+1;\phi)}{{}_2F_1(-\gamma,\beta,\alpha+\beta;\phi)}\right)^2.$$

From Theorem 5 we can see that the moments of the HyperBeta distribution also involve the hypergeometric function, just as with the density function. This means that the properties of the moments are relatively complex and will involve the same computational challenges as computation of the density function. (We discuss computation in a later section.)

Substitution of the appropriate posterior parameters in the contaminated binomial model (from Tables 2-3) gives posterior moments for those models and it is also possible to establish that these models lead to posterior convergence of the standard form in Bayesian analysis for an IID model. Specifically, taking $\theta_0$ and $\phi_0$ as the true values[4] of the parameters $\theta$ and $\phi$, the model leads to convergence towards posterior concentration around the parameter values with the smallest Kullback-Leibler (KL) divergence from the true model (see Berk 1966, Bunke and Milhaud 1998, Kleijn and van der Vaart 2006, Lee and MacEachern 2011, Ramamoorthi et al 2015). In Theorem 6 we show the KL divergence for the two versions of the contaminated model and the value of $\theta$ that is closest to the true model (even under mis-specification). In both cases, if $\phi = \phi_0$ (i.e., if the model is correctly specified) then we get zero KL divergence when $\theta = \theta_0$. This is a comforting property which shows that correct specification of the contaminated binomial model leads to sensible properties for the KL divergence.

[4] We remind the reader that the parameter $\phi$ has an assumed known value in the contaminated binomial model. By specifying a true value $\phi_0$ we allow for the possibility that the model has been mis-specified (i.e., $\phi \neq \phi_0$).

**THEOREM 6A (KL divergence for contaminated binomial model):** For the contaminated binomial model in Table 2 the Kullback-Leibler divergence between the distributions with parameters $(\theta_0, \phi_0)$ and $(\theta, \phi)$ is:

$$\mathrm{KL}(\theta_0, \phi_0 | \theta, \phi) = \sum_{s=0}^{n} \mathrm{Bin}(s|n, \theta_0\phi_0) \begin{bmatrix} s \log\left(\dfrac{\theta_0\phi_0}{\theta\phi}\right) \\ +(n-s)\log\left(\dfrac{1-\theta_0\phi_0}{1-\theta\phi}\right) \end{bmatrix}.$$

The KL divergence is convex in $\theta$ and is minimised at the point $\theta = [\theta_0\phi_0/\phi]_{0,1}$,[5] yielding zero divergence when $\theta = \theta_0\phi_0/\phi$.

**THEOREM 6B (KL divergence for contaminated binomial model):** For the contaminated binomial model in Table 3 the Kullback-Leibler divergence between the distributions with parameters $(\theta_0, \phi_0)$ and $(\theta, \phi)$ is:

$$\mathrm{KL}(\theta_0, \phi_0 | \theta, \phi) = \sum_{s=0}^{n} \mathrm{Bin}(s|n, (1-\theta_0)\phi_0) \begin{bmatrix} s \log\left(\dfrac{(1-\theta_0)\phi_0}{(1-\theta)\phi}\right) \\ +(n-s)\log\left(\dfrac{1-\phi_0+\theta_0\phi_0}{1-\phi+\theta\phi}\right) \end{bmatrix}.$$

The KL divergence is convex in $\theta$ and is minimised at the point $\theta = [(\phi - \phi_0 + \theta_0\phi_0)/\phi]_{0,1}$,[5] yielding zero divergence when $\theta = (\phi - \phi_0 + \theta_0\phi_0)/\phi$.

Using our results for the KL divergence we can proceed to look at posterior convergence and the convergence of the posterior moments. We again take $\theta_0$ and $\phi_0$ to be the true values of the parameters $\theta$ and $\phi$. In the contaminated binomial model set out in Table 2 above (using the left-HyperBeta distribution) we have the posterior density:

$$\pi(\theta|\boldsymbol{x}_n) = \mathrm{LHyperBeta}(\theta|\alpha + \dot{x}_n, \beta, \gamma + n - \dot{x}_n, \phi) \qquad \dot{x}_n \sim \mathrm{Bin}(n, \theta_0\phi_0),$$

and in the alternative contaminated binomial model set out in Table 3 above (using the right-HyperBeta distribution) we have the posterior density:

[5] The notation $[\,\cdot\,]_{0,1}$ is used here to refer to the "clamped" value of the argument within the range $[0,1]$ (a variation on Macaulay bracket notation) and is formally defined by $[r]_{0,1} \equiv \max(0, \min(r, 1))$ for all $r \in \mathbb{R}$.

$$\pi(\theta|\boldsymbol{x}_n) = \text{RHyperBeta}(\theta|\alpha, \beta + \dot{x}_n, \gamma + n - \dot{x}_n, \phi) \qquad \dot{x}_n \sim \text{Bin}(n, (1-\theta_0)\phi_0)\,.$$

In either case, it is useful to consider the asymptotic behaviour of the posterior and its moments under specification of the true parameter value. Taking $\theta_0$ and $\phi_0$ to be the true values of the parameters, we will show that the posterior distribution for $\theta$ converges towards a point-mass distribution on the value:[6]

$$\theta_* = \begin{cases} [\theta_0\phi_0/\phi]_{0,1} & \text{if } \kappa = 0, \\ [(\phi - \phi_0 + \theta_0\phi_0)/\phi]_{0,1} & \text{if } \kappa = 1. \end{cases}$$

To do this, let $\mathcal{U}_*$ be a neighbourhood of the value $\theta_*$. Under the model our calculation of the posterior probability that the parameter is in the stipulated neighbourhood is:

$$\pi(\mathcal{U}_*|\boldsymbol{x}_n) \equiv \int\limits_{\mathcal{U}_*} \pi(\theta|\boldsymbol{x}_n)\, d\theta\,.$$

Assuming that the specification of the sampling distribution is accurate (so the sample total has the sampling distribution we have stipulated) we can establish that this posterior probability that the parameter is in this neighbourhood converges to one almost surely. We establish this formally for both versions of the contaminated binomial model in Theorem 7.

**THEOREM 7 (Posterior consistency in contaminated binomial model):** In the contaminated binomial model, suppose we take $\theta_0$ and $\phi_0$ to be the true values of the parameters $\theta$ and $\phi$ and let $\mathcal{U}_*$ be a neighbourhood of the value $\theta_*$. Then we have:

$$\mathbb{P}\left(\lim_{n\to\infty} \pi(\mathcal{U}_*|\boldsymbol{x}_n) = 1 \middle| \theta = \theta_0\right) = 1.$$

Theorem 7 is proved as an application of general posterior consistency results for Bayesian analysis set out in Schwarz (1965) and expanded in LeCam (1973) and van der Vaart (1998). For formal purposes this is all that is required, but the proof sidesteps the particular form of the moments for the HyperBeta distribution shown in Theorem 5 above. To augment these results, the reader may find it edifying to see a simple heuristic demonstration of why the form of the moments in Theorem 5 leads to the asserted posterior converge — in particular, why these moment forms give asymptotic convergence $\mathbb{E}(\theta|\boldsymbol{x}_n) \to \theta_*$ and $\mathbb{V}(\theta|\boldsymbol{x}_n) \to 0$ as part of the stronger convergence result. We give a heuristic demonstration of this result for the first contaminated binomial model (using the left-HyperBeta distribution).

[6] Note again the use of the "clamped" bracket described in footnote 3.

Consider the first contaminated binomial model in Table 2 and take $\theta_0$ and $\phi_0$ to be the true values of the parameters $\theta$ and $\phi$. (We note that we will still assume that the model uses $\phi$ as the assumed push proportion, so the model is mis-specified if $\phi \neq \phi_0$.) Using Theorem 6 and given the observed sample vector $\boldsymbol{x}_n$ the posterior raw moments for the model are:

$$\mathbb{E}(\theta^k|\boldsymbol{x}_n) = \frac{\int_0^1 \theta^{\alpha+k+\dot{x}_n-1}(1-\theta)^{\beta-1}(1-\theta\phi)^{\gamma+n-\dot{x}_n} d\theta}{\int_0^1 \theta^{\alpha+\dot{x}_n-1}(1-\theta)^{\beta-1}(1-\theta\phi)^{\gamma+n-\dot{x}_n} d\theta}.$$

To facilitate heuristic analysis, suppose we define the function $H_{n,k}: [0,1] \to \mathbb{R}$ by:

$$H_{n,k}(\theta) = \theta^{\alpha+k+n\theta_0\phi_0-1}(1-\theta\phi)^{\gamma+n(1-\theta_0\phi_0)} \qquad 0 \leq \theta \leq 1.$$

In this model we have $\dot{x}_n \sim \text{Bin}(n, \theta_0\phi_0)$ so we have convergence $\dot{x}_n/n \to \theta_0\phi_0$ as $n \to \infty$ (in a sense that we leave unspecified for heuristic purposes). Suppose we consider the case of large $n$ and take the asymptotic equivalence $\dot{x}_n \cong n\theta_0\phi_0$. If we substitute this into the formula for the posterior moments we get the asymptotic equivalence:

$$\begin{aligned}\mathbb{E}(\theta^k|\boldsymbol{x}_n) &\cong \frac{\int_0^1 \theta^{\alpha+k+n\theta_0\phi_0-1}(1-\theta)^{\beta-1}(1-\theta\phi)^{\gamma+n(1-\theta_0\phi_0)} d\theta}{\int_0^1 \theta^{\alpha+n\theta_0\phi_0-1}(1-\theta)^{\beta-1}(1-\theta\phi)^{\gamma+n(1-\theta_0\phi_0)} d\theta} \\ &= \frac{\int_0^1 H_{n,k}(\theta)(1-\theta)^{\beta-1} d\theta}{\int_0^1 H_{n,0}(\theta)(1-\theta)^{\beta-1} d\theta}.\end{aligned}$$

With a bit of calculus, it can easily be shown that the function $H_{n,k}$ is strictly concave with a unique maximising point and resulting maxima given respectively by:

$$m_{n,k} \equiv \underset{0\leq\theta\leq1}{\arg\max}\, H_{n,k}(\theta) = \frac{1}{\phi} \cdot \frac{\alpha+k+n\theta_0\phi_0-1}{\alpha+\gamma+n+k-1},$$

$$M_{n,k} \equiv \max_{0\leq\theta\leq1} H_{n,k}(\theta) = \frac{(\alpha+k+n\theta_0\phi_0-1)^{\alpha+k+n\theta_0\phi_0-1}(\gamma+n(1-\theta_0\phi_0))^{\gamma+n(1-\theta_0\phi_0)}}{\phi^{\alpha+k+n\theta_0\phi_0-1}(\alpha+\gamma+n+k-1)^{\alpha+k+\gamma+n-1}}.$$

As $n \to \infty$ we have $m_{n,k} \to \theta_*$ and $H_{n,k}(\theta)/M_{n,k} \to \mathbb{I}(\theta = \theta_*)$ which means that the function $H_{n,k}$ vanishes asymptotically relative to its maxima at all points except at the limit of the maximising point. Since this function is in the integrand of the numerator and denominator of the posterior moment, and it is continuous, as $n \to \infty$ we obtain the asymptotic equivalence:

$$\begin{aligned}\mathbb{E}(\theta^k|\boldsymbol{x}_n) &\cong \frac{\int_0^1 H_{n,k}(\theta)(1-\theta)^{\beta-1} d\theta}{\int_0^1 H_{n,0}(\theta)(1-\theta)^{\beta-1} d\theta} \\ &= \frac{M_{n,k}}{M_{n,0}} \cdot \frac{\int_0^1 (H_{n,k}(\theta)/M_{n,k})(1-\theta)^{\beta-1} d\theta}{\int_0^1 (H_{n,0}(\theta)/M_{n,0})(1-\theta)^{\beta-1} d\theta}\end{aligned}$$

$$
\begin{aligned}
&\cong \frac{M_{n,k}}{M_{n,0}} \cdot \frac{(1-\theta_*)^{\beta-1}}{(1-\theta_*)^{\beta-1}} \\
&= \frac{1}{\phi^k} \cdot \frac{(\alpha+k+n\theta_0\phi_0-1)^{\alpha+k+n\theta_0\phi_0-1}}{(\alpha+n\theta_0\phi_0-1)^{\alpha+n\theta_0\phi_0-1}} \cdot \frac{(\alpha+\gamma+n-1)^{\alpha+\gamma+n-1}}{(\alpha+\gamma+n+k-1)^{\alpha+k+\gamma+n-1}} \\
&= \frac{1}{\phi^k} \cdot \left(1+\frac{k}{\alpha+n\theta_0\phi_0-1}\right)^{\alpha+n\theta_0\phi_0-1} \cdot \left(1-\frac{k}{\alpha+\gamma+n+k-1}\right)^{\alpha+\gamma+n-1} \\
&\qquad \times \left(\frac{\alpha+k+n\theta_0\phi_0-1}{\alpha+\gamma+n+k-1}\right)^k \\
&\to \frac{1}{\phi^k} \cdot \exp(k) \cdot \exp(-k) \times (\theta_0\phi_0)^k \\
&= \left(\theta_0 \cdot \frac{\phi_0}{\phi}\right)^k \\
&= \theta_*^k.
\end{aligned}
$$

From this general moment convergence result we then have:

$$
\mathbb{E}(\theta|\boldsymbol{x}_n) \to \theta_* \qquad \mathbb{V}(\theta|\boldsymbol{x}_n) \to 0.
$$

This shows that the posterior converges in mean-square (and therefore also in probability) to a point-mass distribution on $\theta_*$.[7] In the special case where $\phi = \phi_0$ (i.e., where the model is correctly specified) we have the posterior converges to a point-mass distribution on $\theta_0$ so we have posterior consistency for the parameter of interest. (The corresponding heuristic argument for posterior convergence in the alternative contaminated binomial model using the right-HyperBeta distribution is analogous; we omit it for brevity.)

From Theorem 7 we see that if the contaminated binomial model is correctly specified (i.e., if the stipulated parameter $\phi$ is equal to its true value) then there is posterior convergence to the true value of the parameter $\theta$ (i.e., the model displays posterior consistency). The theorem also shows the consequences of mis-specification, which is posterior inconsistency in which the posterior converges to a point mass on the broader value $\theta_*$. It should be unsurprising that the contaminated binomial model hinges on correct specification of the probability $\phi$, and the practical lesson is that the model should generally only be used in cases where $\phi$ is a control parameter (i.e., a parameter under the control of the analyst) so that it is known.

[7] Of course, the convergence result in Theorem 8 is stronger than this, but hopefully the above heuristic argument augments the intuition of this result.

## 4. Other properties of the directional HyperBeta distributions

In this section we examine some further properties of the HyperBeta distribution relating to other moment quantities. In particular, we will examine the entropy of the distribution and the score function for the density, which each use functions formed by taking expected values of relevant logarithmic terms in the density. To facilitate this analysis, consider a random variable $X_L \sim \text{LHyperBeta}(x|\alpha,\beta,\gamma,\phi)$ and define the expected value functions:

$$L_1(\alpha,\beta,\gamma,\phi) \equiv \mathbb{E}(\log(X_L)) = \int_0^1 \log(x)\,\text{LHyperBeta}(x|\alpha,\beta,\gamma,\phi)\,dx\,,$$

$$L_2(\alpha,\beta,\gamma,\phi) \equiv \mathbb{E}(\log(1-X_L)) = \int_0^1 \log(1-x)\,\text{LHyperBeta}(x|\alpha,\beta,\gamma,\phi)\,dx\,,$$

$$L_3(\alpha,\beta,\gamma,\phi) \equiv \mathbb{E}(\log(1-X_L\phi)) = \int_0^1 \log(1-x\phi)\,\text{LHyperBeta}(x|\alpha,\beta,\gamma,\phi)\,dx\,,$$

and let $X_R \sim \text{RHyperBeta}(x|\alpha,\beta,\gamma,\phi)$ and define the expected value functions:

$$R_1(\alpha,\beta,\gamma,\phi) \equiv \mathbb{E}(\log(X_R)) = \int_0^1 \log(x)\,\text{RHyperBeta}(x|\alpha,\beta,\gamma,\phi)\,dx\,,$$

$$R_2(\alpha,\beta,\gamma,\phi) \equiv \mathbb{E}(\log(1-X_R)) = \int_0^1 \log(1-x)\,\text{RHyperBeta}(x|\alpha,\beta,\gamma,\phi)\,dx\,,$$

$$R_3(\alpha,\beta,\gamma,\phi) \equiv \mathbb{E}(\log(1-\phi+X_R\phi)) = \int_0^1 \log(1-\phi+x\phi)\,\text{RHyperBeta}(x|\alpha,\beta,\gamma,\phi)\,dx\,.$$

These functions allow us to write the Shannon entropy of the HyperBeta distribution, which is shown in Theorem 8. They also allow us to write the partial derivatives of the log-density function which means that they appear in the score function when taking IID observations from the distribution, which is shown in Theorem 9.

**THEOREM 8A (Left-HyperBeta entropy):** Taking $X \sim \text{LHyperBeta}(\alpha,\beta,\gamma,\phi)$ gives:

$$\text{H}(X) = \mathbb{E}(-\log X) = -L_1(\alpha,\beta,\gamma,\phi).$$

**THEOREM 8B (Right-HyperBeta entropy):** Taking $X \sim \text{RHyperBeta}(\alpha,\beta,\gamma,\phi)$ gives:

$$\text{H}(X) = \mathbb{E}(-\log X) = -R_1(\alpha,\beta,\gamma,\phi).$$

**THEOREM 9A (Log-likelihood and score function for the left-HyperBeta distribution):** If we take the values $X_1, \dots, X_n \sim \text{IID LHyperBeta}(\alpha, \beta, \gamma, \phi)$ and observe the corresponding sample vector $\boldsymbol{x}_n = (x_1, \dots, x_n)$ then the resulting log-likelihood function is:

$$\ell_{\boldsymbol{x}_n}(\alpha, \beta, \gamma, \phi, 0) = (\alpha - 1)\sum_{i=1}^{n} \log(x_i) + (\beta - 1)\sum_{i=1}^{n} \log(1 - x_i) + \gamma \sum_{i=1}^{n} \log(1 - x_i\phi) - \sum_{i=1}^{n} \log\left(\int_0^1 x_i^{\alpha-1}(1 - x_i)^{\beta-1}(1 - x_i\phi)^{\gamma} dx\right),$$

and the resulting partial derivatives in the score function are:

$$\frac{\partial \ell_{\boldsymbol{x}_n}}{\partial \alpha}(\alpha, \beta, \gamma, \phi) = \sum_{i=1}^{n} \log(x_i) - nL_1(\alpha, \beta, \gamma, \phi),$$

$$\frac{\partial \ell_{\boldsymbol{x}_n}}{\partial \beta}(\alpha, \beta, \gamma, \phi) = \sum_{i=1}^{n} \log(1 - x_i) - nL_2(\alpha, \beta, \gamma, \phi),$$

$$\frac{\partial \ell_{\boldsymbol{x}_n}}{\partial \gamma}(\alpha, \beta, \gamma, \phi) = \sum_{i=1}^{n} \log(1 - x_i\phi) - nL_3(\alpha, \beta, \gamma, \phi),$$

$$\frac{\partial \ell_{\boldsymbol{x}_n}}{\partial \phi}(\alpha, \beta, \gamma, \phi) = -\sum_{i=1}^{n} \frac{\gamma x_i}{1 - x_i\phi} + n\gamma \cdot \frac{\text{B}(\alpha, \beta)\ {}_2F_1(-\gamma, \alpha, \alpha + \beta; \phi)}{\text{B}(\alpha + 1, \beta)\ {}_2F_1(-\gamma - 1, \alpha + 1, \alpha + \beta + 1; \phi)}.$$

**THEOREM 9B (Log-likelihood and score function for the right-HyperBeta distribution):** If we take the values $X_1, \dots, X_n \sim \text{IID RHyperBeta}(\alpha, \beta, \gamma, \phi)$ and observe the corresponding sample vector $\boldsymbol{x}_n = (x_1, \dots, x_n)$ then the resulting log-likelihood function is:

$$\ell_{\boldsymbol{x}_n}(\alpha, \beta, \gamma, \phi, 1) = (\alpha - 1)\sum_{i=1}^{n} \log(x_i) + (\beta - 1)\sum_{i=1}^{n} \log(1 - x_i) + \gamma \sum_{i=1}^{n} \log(1 - \phi + x_i\phi) - \sum_{i=1}^{n} \log\left(\int_0^1 x_i^{\alpha-1}(1 - x_i)^{\beta-1}(1 - \phi + x_i\phi)^{\gamma} dx\right),$$

and the resulting partial derivatives in the score function are:

$$\frac{\partial \ell_{\boldsymbol{x}_n}}{\partial \alpha}(\alpha, \beta, \gamma, \phi) = \sum_{i=1}^{n} \log(x_i) - nR_1(\alpha, \beta, \gamma, \phi),$$

$$\frac{\partial \ell_{\boldsymbol{x}_n}}{\partial \beta}(\alpha, \beta, \gamma, \phi) = \sum_{i=1}^{n} \log(1 - x_i) - nR_2(\alpha, \beta, \gamma, \phi),$$

$$\frac{\partial \ell_{\boldsymbol{x}_n}}{\partial \gamma}(\alpha, \beta, \gamma, \phi) = \sum_{i=1}^{n} \log(1 - \phi + x_i\phi) - nR_3(\alpha, \beta, \gamma, \phi),$$

$$\frac{\partial \ell_{x_n}}{\partial \phi}(\alpha,\beta,\gamma,\phi) = -\sum_{i=1}^{n} \frac{\gamma(1-x_i)}{1-\phi+x_i\phi} + n\gamma \cdot \frac{\mathrm{B}(\alpha,\beta)\ \ {}_2F_1(-\gamma,\alpha,\alpha+\beta;\phi)}{\mathrm{B}(\alpha,\beta+1)\ \ {}_2F_1(-\gamma-1,\alpha,\alpha+\beta+1;\phi)}.$$

Our analysis in this paper has focused on Bayesian modelling, but for the sake of completeness it is useful to note some results from classical estimation that come out of the above forms.[8] The above results can be used to obtain the maximum likelihood estimator for the parameters using numerical methods (e.g., via gradient descent). It is notable that if we condition on $\phi$ then the optima for the other parameters are effectively method-of-moments estimators. We will illustrate this for the left-HyperBeta distribution but the right case is analogous. For the left-HyperBeta distribution, taking any fixed value for $\phi$, the other distribution parameters can be maximised (conditional on $\phi$) by setting the first three partial derivatives in the score function to zero, which is equivalent to solving the method-of-moments equations:

$$\frac{1}{n}\sum_{i=1}^{n} \log(x_i) = \mathbb{E}(\log(X_L)),$$

$$\frac{1}{n}\sum_{i=1}^{n} \log(1-x_i) = \mathbb{E}(\log(1-X_L)),$$

$$\frac{1}{n}\sum_{i=1}^{n} \log(1-x_i\phi) = \mathbb{E}(\log(1-X_L\phi)).$$

This means that the conditional optima $(\hat{\alpha},\hat{\beta},\hat{\gamma})$ (given the value of $\phi$) are method-of-moments estimators (which can be computed numerically through gradient methods). (If the parameter $\phi$ is taken to be unknown then it is possible to find its MLE from the remaining score equation, but this is not necessary in our analysis, since we assume the parameter to be known.)

[8] Since the HyperBeta distribution arises naturally from a Bayesian analysis of the contaminated binomial model, it is uncommon to use the MLE to estimate the parameters. Instead, it is usual to estimate via quantities such as the MAP estimator, which corresponds to the posterior mode of the pushed beta distribution. In any case, the MLE is here for convenience.

## 5. Computing with the directional Hyperbeta distributions

Many of the functions describing the HyperBeta distribution cannot be written in closed form and they require numerical methods to compute effectively. In particular, the density function, cumulative distribution function and moment functions all require integration of the density kernel and the quantile function and pseudo-random generation function also requires inversion of the resulting values. The kernel function is non-negative and sometimes gives very small values, so it is useful to compute it in log-space via the following log-integral functions (for the left and right-pushed distributions respectively):

$$H_0(r|\alpha,\beta,\gamma,\phi) \equiv \log\left(\int_0^r x^{\alpha-1}(1-x)^{\beta-1}(1-x\phi)^{\gamma}dx\right),$$

$$H_1(r|\alpha,\beta,\gamma,\phi) \equiv \log\left(\int_0^r x^{\alpha-1}(1-x)^{\beta-1}(1-\phi+x\phi)^{\gamma}dx\right).$$

The full integrals in these log-integral functions (with $r = 1$) are closely related to the Gaussian hypergeometric function ${}_2F_1$ and the beta function B, via the equations:

$$H_0(1|\alpha,\beta,\gamma,\phi) = \log \mathrm{B}(\alpha,\beta) + \log\ {}_2F_1(-\gamma,\alpha,\alpha+\beta;\phi)\,,$$

$$H_1(1|\alpha,\beta,\gamma,\phi) = \log \mathrm{B}(\alpha,\beta) + \log\ {}_2F_1(-\gamma,\beta,\alpha+\beta;\phi)\,.$$

The partial integrals have a corresponding relationship to the partial Gaussian hypergeometric function and the beta function.

As a starting point for computation we use the default **`integrate`** function in **`R`**, which uses globally adaptive interval subdivision using Gauss–Kronrod quadrature (see Notaris 2016).[9] This function is based on the Quadpack routines **`dqags`** and **`dqagi`** in Piessens *et al* (1989). The **`integrate`** function does a serviceable job of computing the above log-integrals for most inputs of interest, but computational underflow causes it to fail to compute a positive value (i.e., it computes the value at zero) for the integrals when the shape and/or intensity parameters are so large that the kernel function too peaked.[10]

[9] See documentation at https://stat.ethz.ch/R-manual/R-devel/library/stats/html/integrate.html

[10] Investigations by the present author show that the function fails to compute the integral accurately when the inputs for the shape and intensity parameters are in the thousands or higher — in such cases the integration yields a computed value of zero and so the log-integral is given as negative infinity (when the true value is higher than this). The documentation for the **`integrate`** function warns specifically about this deficiency, noting that: "If the function is approximately constant (in particular, zero) over nearly all its range it is possible that the result and error estimate may be seriously wrong."

To get around this problem, we use a midpoint quadrature method based on sample points taken from the quantile function for the beta distribution, which concentrate near to the peak of the kernel function, reducing arithmetic underflow in the integral approximation. To do this we first define the function:

$$S(x|\phi) \equiv \log(1 - x\phi) - \log(1 - x).$$

Our method is based on the observation that the integrals of interest can be written as:

$$
\begin{aligned}
\exp(H_0(r|\alpha,\beta,\gamma,\phi)) &= \int_0^r x^{\alpha-1}(1-x)^{\beta-1}(1-x\phi)^{\gamma} dx \\
&= \int_0^r x^{\alpha-1}(1-x)^{\beta+\gamma-1}\left(\frac{1-x\phi}{1-x}\right)^{\gamma} dx \\
&= \frac{\Gamma(\alpha)\Gamma(\beta+\gamma)}{\Gamma(\alpha+\beta+\gamma)}\int_0^r \text{Beta}(x|\alpha,\beta+\gamma)\left(\frac{1-x\phi}{1-x}\right)^{\gamma} dx \\
&= \frac{\Gamma(\alpha)\Gamma(\beta+\gamma)}{\Gamma(\alpha+\beta+\gamma)}\int_0^r \exp(\gamma S(x|\phi)) \cdot \text{Beta}(x|\alpha,\beta+\gamma)\, dx\,, \\
&= \frac{\Gamma(\alpha)\Gamma(\beta+\gamma)}{\Gamma(\alpha+\beta+\gamma)}\int_0^1 \mathbb{I}(x \leq r) \cdot \exp(\gamma S(x|\phi)) \cdot \text{Beta}(x|\alpha,\beta+\gamma)\, dx\,,
\end{aligned}
$$

$$
\begin{aligned}
\exp(H_1(r|\alpha,\beta,\gamma,\phi)) &= \int_0^r x^{\alpha-1}(1-x)^{\beta-1}(1-\phi+x\phi)^{\gamma} dx \\
&= \int_0^r x^{\alpha+\gamma-1}(1-x)^{\beta-1}\left(\frac{1-(1-x)\phi}{1-(1-x)}\right)^{\gamma} dx \\
&= \frac{\Gamma(\alpha+\gamma)\Gamma(\beta)}{\Gamma(\alpha+\beta+\gamma)}\int_0^r \text{Beta}(1-x|\alpha+\gamma,\beta)\left(\frac{1-(1-x)\phi}{1-(1-x)}\right)^{\gamma} dx \\
&= \frac{\Gamma(\alpha+\gamma)\Gamma(\beta)}{\Gamma(\alpha+\beta+\gamma)}\int_{1-r}^1 \text{Beta}(x|\beta,\alpha+\gamma)\left(\frac{1-x\phi}{1-x}\right)^{\gamma} dx \\
&= \frac{\Gamma(\alpha)\Gamma(\beta+\gamma)}{\Gamma(\alpha+\beta+\gamma)}\int_{1-r}^1 \exp(\gamma S(x|\phi)) \cdot \text{Beta}(x|\beta,\alpha+\gamma)\, dx \\
&= \frac{\Gamma(\alpha)\Gamma(\beta+\gamma)}{\Gamma(\alpha+\beta+\gamma)}\int_0^1 \mathbb{I}(x \geq 1-r) \cdot \exp(\gamma S(x|\phi)) \cdot \text{Beta}(x|\beta,\alpha+\gamma)\, dx\,,
\end{aligned}
$$

For the left-HyperBeta integral, we generate the values $q_1, \dots, q_M \sim \text{Beta}(\alpha, \beta+\gamma)$ and write the ordered values as $q_{(1)} < \cdots < q_{(M)}$. To represent the indicator function in the integral we form the corresponding log-weight values:

$$w_i \equiv \begin{cases} 0 & \text{if } q_{(i)} \leq r, \\ \log(r - q_{(i-1)}) - \log(q_{(i)} - q_{(i-1)}) & \text{if } q_{(i-1)} \leq r < q_{(i)}, \\ -\infty & \text{if } q_{(i-1)} > r. \end{cases}$$

We then obtain the integral estimator:

$$\exp(\widehat{H}_0(r|\alpha,\beta,\gamma,\phi)) \approx \frac{\Gamma(\alpha)\Gamma(\beta+\gamma)}{\Gamma(\alpha+\beta+\gamma)} \cdot \frac{1}{M} \sum_{i=1}^{M} \exp(w_i + \gamma S(q_{(i)}|\phi)),$$

which gives the corresponding log-space equation:

$$\begin{aligned} \widehat{H}_0(r|\alpha,\beta,\gamma,\phi) \approx & \operatorname*{logsumexp}_i(w_i + \gamma S(q_{(i)}|\phi)) - \log(M) \\ & + \log\Gamma(\alpha) + \log\Gamma(\beta+\gamma) - \log\Gamma(\alpha+\beta+\gamma). \end{aligned}$$

For the right- HyperBeta integral, we generate the values $q_1, \dots, q_M \sim \text{Beta}(\beta, \alpha+\gamma)$ and write the ordered values as $q_{(1)} < \cdots < q_{(M)}$. To represent the indicator function in the integral we form the corresponding log-weight values:

$$w_i \equiv \begin{cases} -\infty & \text{if } q_{(i)} < 1-r, \\ \log(q_{(i)} - 1 + r) - \log(q_{(i)} - q_{(i-1)}) & \text{if } q_{(i-1)} < 1-r \leq q_{(i)}, \\ 0 & \text{if } q_{(i-1)} \geq 1-r. \end{cases}$$

We then obtain the integral estimator:

$$\exp(\widehat{H}_1(r|\alpha,\beta,\gamma,\phi)) \approx \frac{\Gamma(\alpha+\gamma)\Gamma(\beta)}{\Gamma(\alpha+\beta+\gamma)} \cdot \frac{1}{M} \sum_{i=1}^{M} \exp(w_i + \gamma S(q_{(i)}|\phi)),$$

which gives the corresponding log-space equation:

$$\begin{aligned} \widehat{H}_1(r|\alpha,\beta,\gamma,\phi) \approx & \operatorname*{logsumexp}_i(w_i + \gamma S(q_{(i)}|\phi)) - \log(M) \\ & + \log\Gamma(\alpha+\gamma) + \log\Gamma(\beta) - \log\Gamma(\alpha+\beta+\gamma). \end{aligned}$$

**REMARK:** The log-weight terms $w_i$ incorporate a weighting into the calculation representing the indicator functions for the range of the integral. It is possible to proceed in an alternative way by generating the importance sample over the range $[0, r]$ instead of adding a weighting term. We have chosen to generate the approximation using a set of importance sample values over the whole unit interval so that the log-kernel function can be vectorised for an input vector of values for $r$ and the resulting approximations will be linear piecewise consistent. □

In either of the above cases we can form a deterministic approximation to the integral by taking the $q_i$ to be evenly spaced quantiles with corresponding probabilities $p_i = (2i-1)/2M$. The quantiles tend to clump in the peaks of the density function which reduces arithmetic underflow in the computation of the integral. Investigations by the present author show that this integral approximation method performs well when the input parameters are large and the resulting density function is highly peaked (the case when naïve integral methods that sample uniformly over the unit interval fail). Computation of the log-kernel function above is implemented in the **`scale.hyperbeta`** function in **R** (O'Neill 2025). This function uses the **`integrate`** function for base computation but switches to the method shown above if the former method gives a zero value for the integral (which is an incorrect value). In the latter case the function uses $M = 10^6$ quantile values by default.

Once we have a computational facility that can compute the log-kernel function effectively we can then write the various probability functions for the pushed beta distribution. The density function, cumulative distribution function and raw moments are given on the log-scale as:

$$\begin{aligned}
\log \text{HyperBeta}(x|\alpha,\beta,\gamma,\phi,\kappa) &= (\alpha-1)\log(x) + (\beta-1)\log(1-x) - H_\kappa(1|\alpha,\beta,\gamma,\phi) \\
&\quad + (1-\kappa)\gamma\log(1-x\phi) + \kappa\gamma\log(1-\phi+x\phi)\,, \\
\log \text{CDF}_{\text{HyperBeta}}(x|\alpha,\beta,\gamma,\phi,\kappa) &= H_\kappa(x|\alpha,\beta,\gamma,\phi) - H_\kappa(1|\alpha,\beta,\gamma,\phi), \\
\log \mathbb{E}(X^k) &= H_\kappa(1|\alpha+k,\beta,\gamma,\phi) - H_\kappa(1|\alpha,\beta,\gamma,\phi).
\end{aligned}$$

The quantile function can be computed from the cumulative distribution using appropriate root-finding methods and the pseudo-random generation function can then be computed using the quantile function using inverse-transform sampling.

The directional versions of the HyperBeta distribution is implemented in **R** in a set of functions available in a public GitHub repository (O'Neill 2025). These functions use the directional form of the distribution allowing use of the left-HyperBeta and right-HyperBeta variants of the distribution. Table 4 below shows the available functions and their inputs. The available functions are the density function, cumulative distribution function, quantile function, random generation function, HDR function and moment function. Each of these functions takes the parameters of the distribution given as **`shape1`**, **`shape2`**, **`intensity`** and **`proportion`**. There is also a logical value **`right`** that controls the direction of the push/pull, with the left-HyperBeta distribution used as the default. Other inputs to the functions are standard inputs for the type of probability function being used.

| **TABLE 4: Functions for the HyperBeta distribution** | |
|---|---|
| **Function** | **Inputs** |
| `dhyperbeta` | `x, shape1, shape2, intensity, proportion, right = FALSE, log = FALSE` |
| `phyperbeta` | `x, shape1, shape2, intensity, proportion, right = FALSE, lower.tail = TRUE, log.p = FALSE` |
| `qhyperbeta` | `p, shape1, shape2, intensity, proportion, right = FALSE, lower.tail = TRUE, log.p = FALSE` |
| `rhyperbeta` | `n, shape1, shape2, intensity, proportion, right = FALSE` |
| `HDR.hyperbeta` | `cover.prob, shape1, shape2, intensity, proportion, right = FALSE` |
| `moments.hyperbeta` | `shape1, shape2, intensity, proportion, right = FALSE, include.sd = FALSE` |
| `scale.hyperbeta` | `shape1, shape2, intensity, proportion, right = FALSE, log = TRUE, intvals = 10^6, ...` |

## 6. Applications to contaminated binary sampling

Binary sampling (without contamination) is ubiquitous in research and in statistical problems. It involves solicitation of binary information from participants with subsequent estimation of the prevalence of affirmative outcomes in the population of interest. This can occur in sampling where participants are surveyed on whether or not they possess some characteristic of interest and the goal is to estimate the prevalence of that characteristic in the population.

The binary model is a simple and robust way to estimate the prevalence of a characteristic of interest in a population in the case where there is reason to believe that survey responses are true reflections of the underlying characteristics of participants. However, there are situations in research where the goal is to determine the prevalence of some "sensitive" characteristic where there may be social desirability bias in participant answers to questions that they think may become attributable to them. For example, in certain types of research, participants might be asked whether they have ever been unfaithful to their spouse, have ever committed a serious (and non-detected) crime, have ever had some sexually-transmitted disease, hold some political view, etc. Privacy and confidentiality protections for the survey may go some way towards encouraging truthful answers to these questions, but survey participants may still perceive a danger of "leaked" answers and may be subject to social desirability bias.

In such cases, one potentially valuable experimental protocol is for the researcher to introduce "controlled contamination" of the survey, whereby the response of the participant to a sensitive question is contaminated with an event with known probability (but unknown outcome). This can occur by asking the participant to roll a die, flip a coin, etc., during the survey and to then allow the result of this randomisation device to affect their response to the sensitive question in a stipulated way that "masks" a socially undesirable answer. The basic idea is to ensure that if the participant's answer becomes known (e.g., due to a leak in survey answers) the participant maintains "plausible deniability" against a socially undesirable answer. The contamination allows the participant to answer truthfully (subject to the contamination mechanism) safe in the knowledge that a socially undesirable answer could have come from the contamination mechanism with some substantial probability rather that reflecting their true characteristics. From a statistical perspective, this protocol replaces a source of uncontrolled contamination (social desirability bias) that is hard to model with controlled contamination (a randomisation mechanism with known probabilities) that is easy to model.

**EXAMPLE:** Some social science researchers wish to survey the prevalence of infidelity in a population of interest. To do this, they randomly sample from a sampling frame composed of married people and they administer a survey asking whether or not the subjects have ever been unfaithful to their spouse during their marriage. To elicit reliable information on infidelity (despite the sensitivity), they impose controlled contamination on the question as shown below.

| **SURVEY QUESTION** | |
|---|---|
| *This survey involves "controlled contamination" as an additional protection to your privacy. The controlled contamination means that anyone reading your survey answers will not be able to determine your personal characteristics pertaining to the subject matter of the survey.*<br><br>**Instructions:** Before answering this question, please roll the (six-sided) die in front of you. You do not need to record the outcome you rolled or disclose the outcome to the researcher (i.e., you may keep it a secret only to yourself). If you roll a 1-2, please answer the question below truthfully. If you roll a 3-6, please give the "Yes" answer (even if this is not true). | |
| **Question:** During your marriage (to your present spouse), have you ever been unfaithful? | ☐ Yes (or you rolled 3-6)<br>☐ No |

Observe from the above mechanism that if a participant gives a "Yes" answer to this question, they maintain plausible deniability for infidelity if their survey answer becomes attributable to them. This occurs because there is a significant probability that they simply rolled a 3-6 on the die and therefore gave that answer irrespective of their actual marital conduct. This mechanism therefore constitutes a good protection against social desirability bias (or more direct adverse consequences from an upset wife or husband who sees the answer to the question).

In this situation, the researchers would use the contaminated binomial model in Table 3, leading to a posterior inference from the right-HyperBeta distribution. Suppose we let $X_i$ denote the indicator for a "No" answer for participant $i$. Then the answers to the survey question for the $n$ participants in the survey are contaminated Bernoulli values:[11]

$$X_1, X_2, \ldots, X_n \sim \text{IID Bern}((1-\theta)\phi),$$

where $\theta$ is the true prevalence of infidelity in the population and $\phi = ⅓$ is the probability of a non-contaminated answer. Suppose that the researchers use a uniform prior distribution for $\theta$ and observe $\dot{x}_n = 92$ "No" answers and $n - \dot{x}_n = 248$ "Yes" answers (with $n = 340$). This prior and data would lead to the posterior distribution:

$$\theta|\boldsymbol{x}_n \sim \text{RHyperBeta}(1, 93, 248, ⅓).$$

The posterior expected value and standard deviation for the prevalence parameter are:

$$\mathbb{E}(\theta|\boldsymbol{x}_n) = 0.1856469 \qquad \mathbb{S}(\theta|\boldsymbol{x}_n) = 0.0701663.$$

Thus, based on the posterior expected value, the researchers estimate that the prevalence of infidelity in the population is 18.56%. □

As can be seen from the above example, it is possible to make inferences in cases where there is controlled contamination in binary data. Introducing contamination reduces the information in the observed data since it contaminates the binary indicator of the characteristic of interest. However, this contamination may be worthwhile because it plausibly removes uncontrolled contamination from other sources. It is important to note that there is a trade-off involved in the level of contamination used — the more the contamination mechanism acts in favour of the socially undesirable answer, the more protection is afforded to the participant (since they can more plausibly attribute that answer to the contamination mechanism instead of themselves), but this may also act to reduce the amount of information from the data.

[11] This assumes that the controlled contamination is sufficient to induce truthful answers by survey participants.

## 7. Summary and conclusion

In this paper we have introduced directional variants of the Gaussian hypergeometric beta (HyperBeta) distribution using an alternative parameterisation. This alternative form is useful because it represents the directional push/pull in the distribution in a simple and symmetrical form with a simplified reflection equation for alternating the direction of the push/pull. We have undertaken an extensive examination of the directional HyperBeta distributions and we have shown how they arise naturally under a Bayesian analysis of the contaminated binomial model. The distribution involves an additional density term that pushes/pulls the density to the left or right relative to a beta distribution with the same shape parameters. In this paper we have derived a range of properties of the directional HyperBeta distributions, including shape, modality, moments, and asymptotic properties.

Computation of probabilities and moments of the HyperBeta distribution involves numerical computation of a rather tricky type of integral related to the hypergeometric function. We have examined how to compute the required integral for computation of various quantities of interest using a method that is roughly analogous to importance sampling. This method has been used to create a set of user-friendly functions in **R** that include the probability functions and moment functions for the distribution.

The directional variance of the HyperBeta distribution can be used as conjugate priors in the contaminated binomial model, leading to simple posterior distributions for inference. This model is posterior consistent when the contamination parameter is set correctly (i.e., equal to its true value) but is posterior inconsistent when the contamination parameter is set incorrectly (i.e., not equal to its true value). For this reason, we recommend that contaminated binomial sampling should only be used when the contamination is under the control of the experimenter, such that it is known.

| Statements and Declarations | |
|---|---|
| **Funding** | No funding or other support was received for this research |
| **Conflicts of interest** | The authors have no relevant financial or non-financial interests to disclose |
| **Ethical approval and consent to participate** | Not applicable |
| **Consent for publication** | Not applicable |

## Appendix: Proof of Theorems

**PROOF OF THEOREM 1:** We will prove the first two cases which are for the left-Hyperbeta distribution. In these cases we have $\kappa = 0$ so we consider the kernel function:

$$\text{Kernel}(x|\alpha,\beta,\gamma,\phi) \equiv x^{\alpha-1}(1-x)^{\beta-1}(1-x\phi)^{\gamma} \qquad \text{for } 0 \leq x \leq 1.$$

Define the kernel ratio function using the push-proportions $0 \leq \phi_0 < \phi_1 \leq 1$ by:

$$\text{KR}(x) \equiv \frac{\text{Kernel}(x|\alpha,\beta,\gamma,\phi_0)}{\text{Kernel}(x|\alpha,\beta,\gamma,\phi_1)} = \left(\frac{1-x\phi_0}{1-x\phi_1}\right)^{\gamma}.$$

From the conditions of the theorem we have $\phi_1 - \phi_0 > 0$. Moreover, within the interior of the support $0 < x < 1$ we also have:

$$(1-x\phi_1)^2 > 0 \qquad \frac{1-x\phi_0}{1-x\phi_1} > 1,$$

which then gives the directional result:

$$\frac{d\text{KR}}{dx}(x) = \frac{\gamma(\phi_1-\phi_0)}{(1-x\phi_1)^2}\left(\frac{1-x\phi_0}{1-x\phi_1}\right)^{\gamma-1} > 0 \qquad \gamma > 0 \quad 0 < x < 1,$$

$$\frac{d\text{KR}}{dx}(x) = \frac{\gamma(\phi_1-\phi_0)}{(1-x\phi_1)^2}\left(\frac{1-x\phi_0}{1-x\phi_1}\right)^{\gamma-1} < 0 \qquad \gamma < 0 \quad 0 < x < 1.$$

Now, let us denote the scaling constant $\overline{\text{K}}(\alpha,\beta,\gamma,\phi) \equiv \int_0^1 x^{\alpha-1}(1-x)^{\beta-1}(1-x\phi)^{\gamma}\, dx$ and note that this is strictly positive. We then have:

$$\frac{d}{dx}\frac{\text{LHyperBeta}(x|\alpha,\beta,\gamma,\phi_0)}{\text{LHyperBeta}(x|\alpha,\beta,\gamma,\phi_1)} = \underbrace{\frac{\overline{\text{K}}(\alpha,\beta,\gamma,\phi_1)}{\overline{\text{K}}(\alpha,\beta,\gamma,\phi_0)}}_{+} \cdot \frac{d\text{KR}}{dx}(x) > 0 \qquad \gamma > 0 \quad 0 < x < 1,$$

$$\frac{d}{dx}\frac{\text{LHyperBeta}(x|\alpha,\beta,\gamma,\phi_0)}{\text{LHyperBeta}(x|\alpha,\beta,\gamma,\phi_1)} = \underbrace{\frac{\overline{\text{K}}(\alpha,\beta,\gamma,\phi_1)}{\overline{\text{K}}(\alpha,\beta,\gamma,\phi_0)}}_{+} \cdot \frac{d\text{KR}}{dx}(x) < 0 \qquad \gamma < 0 \quad 0 < x < 1.$$

This establishes that the left-HyperBeta density follows the monotone likelihood property with respect to the push-proportion parameter. The monotone likelihood property implies the first-order stochastic dominance result in the theorem, which establishes the required results. The corresponding proof for the cases where $\kappa = 1$ (the right-Hyperbeta distribution) instead use the term $(1-\phi+x\phi)^{\gamma}$ in the kernel but they are otherwise analogous. ■

**PROOF OF THEOREM 2:** We will prove the first case for the left-Hyperbeta distribution. In this case we have $\kappa = 0$ so we consider the kernel function:

$$\text{Kernel}(x|\alpha,\beta,\gamma,\phi) \equiv x^{\alpha-1}(1-x)^{\beta-1}(1-x\phi)^{\gamma} \qquad \text{for } 0 \leq x \leq 1.$$

Define the kernel ratio function using the push intensity $-\infty < \gamma_0 < \gamma_1 < \infty$ by:

$$\text{KR}(x) \equiv \frac{\text{Kernel}(x|\alpha,\beta,\gamma_0,\phi)}{\text{Kernel}(x|\alpha,\beta,\gamma_1,\phi)} = (1-x\phi)^{-(\gamma_1-\gamma_0)}.$$

From the conditions of the theorem we have $\gamma_1 - \gamma_0 > 0$ and $0 < \phi < 1$. Moreover, within the interior of the support $0 < x < 1$ we also have $1 - x\phi > 0$, which gives the directional result:

$$\frac{d\text{KR}}{dx}(x) = \phi(\gamma_1-\gamma_0)(1-x\phi)^{-(\gamma_1-\gamma_0)-1} > 0 \qquad 0 < x < 1.$$

Now, let us denote the scaling constant $\overline{\text{K}}(\alpha,\beta,\gamma,\phi) \equiv \int_0^1 x^{\alpha-1}(1-x)^{\beta-1}(1-x\phi)^{\gamma}\, dx$ and note that this is strictly positive. We then have:

$$\frac{d}{dx}\frac{\text{LHyperBeta}(x|\alpha,\beta,\gamma_0,\phi)}{\text{LHyperBeta}(x|\alpha,\beta,\gamma_1,\phi)} = \underbrace{\frac{\overline{\text{K}}(\alpha,\beta,\gamma_1,\phi)}{\overline{\text{K}}(\alpha,\beta,\gamma_0,\phi)}}_{+} \cdot \frac{d\text{KR}}{dx}(x) > 0 \qquad 0 < x < 1.$$

This establishes that the left-HyperBeta density follows the monotone likelihood property with respect to the push-intensity parameter. The monotone likelihood property implies the first-order stochastic dominance result in the theorem, which establishes the required results. The corresponding proof for the case where $\kappa = 1$ (the right-Hyperbeta distribution) instead use the term $(1-\phi+x\phi)^{\gamma}$ in the kernel but they are otherwise analogous. ■

**PROOF OF THEOREM 3:** It is simple to establish that the directional HyperBeta distribution obeys the reflection equation:

$$\text{HyperBeta}(x|\alpha,\beta,\gamma,\phi,\kappa) = \text{HyperBeta}(1-x|\beta,\alpha,\gamma,\phi,1-\kappa).$$

Taking $\acute{x} = 1 - x$ and applying this reflection equation gives:

$$\begin{aligned}
\text{DIR}_{\kappa}(x|\alpha,\beta,\gamma,\phi) &= \text{sgn}\frac{d}{dx}\text{HyperBeta}(x|\alpha,\beta,\gamma,\phi,\kappa) \\
&= \text{sgn}\frac{d}{dx}\text{HyperBeta}(1-x|\beta,\alpha,\gamma,\phi,1-\kappa) \\
&= -\,\text{sgn}\frac{d}{d\acute{x}}\text{HyperBeta}(\acute{x}|\beta,\alpha,\gamma,\phi,1-\kappa) \\
&= -\text{DIR}_{1-\kappa}(\acute{x}|\beta,\alpha,\gamma,\phi) \\
&= -\text{DIR}_{1-\kappa}(1-x|\beta,\alpha,\gamma,\phi),
\end{aligned}$$

which was to be shown. ■

**PROOF OF THEOREM 4:** It is established in the body of the paper that for $0 < x < 1$ we have:

$$\text{sgn} \frac{d}{dx} \log \text{HyperBeta}(x|\alpha, \beta, \gamma, \phi, \kappa) = \text{sgn}\, Q_\kappa(x|\alpha, \beta, \gamma, \phi) .$$

Since the logarithm function is a strictly increasing function, the direction of the slope of the HyperBeta density is the same as the direction of the slope of its logarithm we then have:

$$\begin{aligned} \text{DIR}_\kappa(x|\alpha, \beta, \gamma, \phi) &= \text{sgn} \frac{d}{dx} \text{HyperBeta}(x|\alpha, \beta, \gamma, \phi, \kappa) \\ &= \text{sgn} \frac{d}{dx} \log \text{HyperBeta}(x|\alpha, \beta, \gamma, \phi, \kappa) \\ &= \text{sgn}\big[(-1)^\kappa \phi(\alpha + \beta + \gamma - 2)(x - r_{\kappa,0})(x - r_{\kappa,1})\big] \\ &= \text{sgn}[(-1)^\kappa] \times \text{sgn}(\alpha + \beta + \gamma - 2) \times \text{sgn}(x - r_{\kappa,0}) \times \text{sgn}(x - r_{1,0}) . \end{aligned}$$

Substitution of $\kappa = 0$ gives the left case in the theorem. Substitution of $\kappa = 1$ and use of the reflection equation in Theorem 3 then gives the right case in the theorem. ■

**PROOF OF THEOREM 5A:** We first establish the alternative form as:

$$\begin{aligned} \mathbb{E}(X^k) &= \int_0^1 x^k \, \text{HyperBeta}(x|\alpha, \beta, \gamma, \phi) \, dx \\ &= \int_0^1 x^k \cdot \frac{x^{\alpha-1}(1-x)^{\beta-1}(1-x\phi)^\gamma}{\int_0^1 x^{\alpha-1}(1-x)^{\beta-1}(1-x\phi)^\gamma dx} dx \\ &= \frac{\int_0^1 x^{\alpha+k-1}(1-x)^{\beta-1}(1-x\phi)^\gamma dx}{\int_0^1 x^{\alpha-1}(1-x)^{\beta-1}(1-x\phi)^\gamma dx} . \end{aligned}$$

The first form is then established as:

$$\begin{aligned} \mathbb{E}(X^k) &= \int_0^1 x^k \, \text{HyperBeta}(x|\alpha, \beta, \gamma, \phi) \, dx \\ &= \int_0^1 x^k \cdot \frac{x^{\alpha-1}(1-x)^{\beta-1}(1-x\phi)^\gamma}{\text{B}(\alpha, \beta) \ {}_2F_1(-\gamma, \alpha, \alpha+\beta; \phi)} dx \\ &= \frac{\int_0^1 x^{\alpha+k-1}(1-x)^{\beta-1}(1-x\phi)^\gamma dx}{\text{B}(\alpha, \beta) \ {}_2F_1(-\gamma, \alpha, \alpha+\beta; \phi)} \\ &= \frac{\text{B}(\alpha+k, \beta) \ {}_2F_1(-\gamma, \alpha+k, \alpha+\beta+k; \phi)}{\text{B}(\alpha, \beta) \ {}_2F_1(-\gamma, \alpha, \alpha+\beta; \phi)}, \end{aligned}$$

which was to be shown. ■

**PROOF OF THEOREM 5B:** Analogous to proof of Theorem 5A. ■

**PROOF OF THEOREM 6A:** Under the stipulated model we have the sampling density:

$$f(\boldsymbol{x}_n|\theta,\phi) = (\theta\phi)^{\dot{x}_n}(1-\theta\phi)^{n-\dot{x}_n}.$$

This gives the Kullback-Leibler divergence:

$$\begin{aligned}
\text{KL}(\theta_0,\phi_0|\theta,\phi) &= \sum_{\boldsymbol{x}_n} f(\boldsymbol{x}_n|\theta_0,\phi_0)\log\left(\frac{f(\boldsymbol{x}_n|\theta_0,\phi_0)}{f(\boldsymbol{x}_n|\theta,\phi)}\right) \\
&= \sum_{\boldsymbol{x}_n} f(\boldsymbol{x}_n|\theta_0,\phi_0)[\log f(\boldsymbol{x}_n|\theta_0,\phi_0) - \log f(\boldsymbol{x}_n|\theta,\phi)] \\
&= \sum_{\boldsymbol{x}_n} f(\boldsymbol{x}_n|\theta_0,\phi_0)\begin{bmatrix} \dot{x}_n\log(\theta_0\phi_0) + (n-\dot{x}_n)\log(1-\theta_0\phi_0) \\ -\dot{x}_n\log(\theta\phi) - (n-\dot{x}_n)\log(1-\theta\phi) \end{bmatrix} \\
&= \sum_{\boldsymbol{x}_n} f(\boldsymbol{x}_n|\theta_0,\phi_0)\left[\dot{x}_n\log\left(\frac{\theta_0\phi_0}{\theta\phi}\right) + (n-\dot{x}_n)\log\left(\frac{1-\theta_0\phi_0}{1-\theta\phi}\right)\right] \\
&= \sum_{s=0}^{n} \text{Bin}(s|n,\theta_0\phi_0)\left[s\log\left(\frac{\theta_0\phi_0}{\theta\phi}\right) + (n-s)\log\left(\frac{1-\theta_0\phi_0}{1-\theta\phi}\right)\right].
\end{aligned}$$

Now, consider this as a function of $\theta$ written as:

$$R(\theta) \equiv \sum_{s=0}^{n} \text{Bin}(s|n,\theta_0\phi_0)\left[s\log\left(\frac{\theta_0\phi_0}{\theta\phi}\right) + (n-s)\log\left(\frac{1-\theta_0\phi_0}{1-\theta\phi}\right)\right].$$

Taking $S \sim \text{Bin}(n,\theta_0\phi_0)$ the first and second derivatives of this function are:

$$\begin{aligned}
\frac{dR}{d\theta}(\theta) &= \sum_{s=0}^{n} \text{Bin}(s|n,\theta_0\phi_0)\left[-\frac{s}{\theta} + \frac{\phi(n-s)}{1-\theta\phi}\right] \\
&= \mathbb{E}\left[-\frac{S}{\theta} + \frac{\phi(n-S)}{1-\theta\phi}\right] \\
&= \mathbb{E}\left[\frac{-(1-\theta\phi)S + \theta\phi(n-S)}{\theta(1-\theta\phi)}\right] \\
&= \mathbb{E}\left[\frac{n\theta\phi - S}{\theta(1-\theta\phi)}\right] \\
&= \frac{n\theta\phi - \mathbb{E}[S]}{\theta(1-\theta\phi)} \\
&= \frac{n(\theta\phi - \theta_0\phi_0)}{\theta(1-\theta\phi)},
\end{aligned}$$

$$\frac{d^2R}{d\theta^2}(\theta) = \sum_{s=0}^{n} \text{Bin}(s|n,\theta_0\phi_0)\left[\frac{s}{\theta^2} + \frac{\phi^2(n-s)}{(1-\theta\phi)^2}\right]$$

$$= \mathbb{E}\left[\frac{S}{\theta^2} + \frac{\phi^2(n-S)}{(1-\theta\phi)^2}\right]$$
$$= \frac{\mathbb{E}[S]}{\theta^2} + \frac{\phi^2(n-\mathbb{E}[S])}{(1-\theta\phi)^2}$$
$$= \frac{n\theta_0\phi_0}{\theta^2} + \frac{n\phi^2(1-\theta_0\phi_0)}{(1-\theta\phi)^2}$$
$$= n\left[\frac{\theta_0\phi_0}{\theta^2} + \frac{\phi^2(1-\theta_0\phi_0)}{(1-\theta\phi)^2}\right] > 0.$$

This establishes that the KL divergence is strictly convex in $\theta$. Solving the critical equation $0 = dR/d\theta\,(\theta)$ or using the appropriate boundary point (for the monotonic cases) gives the minimising point $\theta = [\theta_0\phi_0/\phi]_{0,1}$. Substituting the case $\theta = \theta_0\phi_0/\phi$ gives:

$$\log\left(\frac{\theta_0\phi_0}{\theta\phi}\right) = \log\left(\frac{1-\theta_0\phi_0}{1-\theta\phi}\right) = 0,$$

which gives zero KL divergence. This establishes all the results to be shown. ■

**PROOF OF THEOREM 6B:** Under the stipulated model we have the sampling density:

$$f(\boldsymbol{x}_n|\theta,\phi) = ((1-\theta)\phi)^{\dot{x}_n}(1-\phi+\theta\phi)^{n-\dot{x}_n}.$$

This gives the Kullback-Leibler divergence:

$$\begin{aligned}
\text{KL}(\theta_0,\phi_0|\theta,\phi) &= \sum_{x_n} f(\boldsymbol{x}_n|\theta_0,\phi_0)\log\left(\frac{f(\boldsymbol{x}_n|\theta_0,\phi_0)}{f(\boldsymbol{x}_n|\theta,\phi)}\right) \\
&= \sum_{x_n} f(\boldsymbol{x}_n|\theta_0,\phi_0)[\log f(\boldsymbol{x}_n|\theta_0,\phi_0) - \log f(\boldsymbol{x}_n|\theta,\phi)]
\end{aligned}$$
$$= \sum_{x_n} f(\boldsymbol{x}_n|\theta_0,\phi_0)\begin{bmatrix} \dot{x}_n\log((1-\theta_0)\phi_0) + (n-\dot{x}_n)\log(1-\phi_0+\theta_0\phi_0) \\ -\dot{x}_n\log((1-\theta)\phi) - (n-\dot{x}_n)\log(1-\phi+\theta\phi) \end{bmatrix}$$
$$= \sum_{x_n} f(\boldsymbol{x}_n|\theta_0,\phi_0)\left[\dot{x}_n\log\left(\frac{(1-\theta_0)\phi_0}{(1-\theta)\phi}\right) + (n-\dot{x}_n)\log\left(\frac{1-\phi_0+\theta_0\phi_0}{1-\phi+\theta\phi}\right)\right]$$
$$= \sum_{s=0}^{n} \text{Bin}(s|n,(1-\theta_0)\phi_0)\left[s\log\left(\frac{(1-\theta_0)\phi_0}{(1-\theta)\phi}\right) + (n-s)\log\left(\frac{1-\phi_0+\theta_0\phi_0}{1-\phi+\theta\phi}\right)\right].$$

Now, consider this as a function of $\theta$ written as:

$$R(\theta) \equiv \sum_{s=0}^{n} \text{Bin}(s|n,(1-\theta_0)\phi_0)\left[s\log\left(\frac{(1-\theta_0)\phi_0}{(1-\theta)\phi}\right) + (n-s)\log\left(\frac{1-\phi_0+\theta_0\phi_0}{1-\phi+\theta\phi}\right)\right].$$

Taking $S \sim \text{Bin}(n,(1-\theta_0)\phi_0)$ the first and second derivatives of this function are:

$$\begin{aligned}
\frac{dR}{d\theta}(\theta) &= \sum_{s=0}^{n} \text{Bin}(s|n,(1-\theta_0)\phi_0)\left[\frac{s}{1-\theta} - \frac{\phi(n-s)}{1-\phi+\theta\phi}\right] \\
&= \mathbb{E}\left[\frac{S}{1-\theta} - \frac{\phi(n-S)}{1-\phi+\theta\phi}\right] \\
&= \mathbb{E}\left[\frac{(1-\phi+\theta\phi)S - (1-\theta)\phi(n-S)}{(1-\theta)(1-\phi+\theta\phi)}\right] \\
&= \mathbb{E}\left[\frac{S - n(1-\theta)\phi}{(1-\theta)(1-\phi+\theta\phi)}\right] \\
&= \frac{\mathbb{E}[S] - n(1-\theta)\phi}{(1-\theta)(1-\phi+\theta\phi)} \\
&= \frac{n[(1-\theta_0)\phi_0 - (1-\theta)\phi]}{(1-\theta)(1-\phi+\theta\phi)},
\end{aligned}$$

$$\begin{aligned}
\frac{d^2R}{d\theta^2}(\theta) &= \sum_{s=0}^{n} \text{Bin}(s|n,\theta_0\phi_0)\left[\frac{s}{(1-\theta)^2} + \frac{\phi^2(n-s)}{(1-\phi+\theta\phi)^2}\right] \\
&= \mathbb{E}\left[\frac{S}{(1-\theta)^2} + \frac{\phi^2(n-S)}{(1-\phi+\theta\phi)^2}\right] \\
&= \frac{\mathbb{E}[S]}{(1-\theta)^2} + \frac{\phi^2(n-\mathbb{E}[S])}{(1-\phi+\theta\phi)^2} \\
&= \frac{n(1-\theta_0)\phi_0}{(1-\theta)^2} + \frac{n\phi^2(1-\phi_0+\theta_0\phi_0)}{(1-\phi+\theta\phi)^2} \\
&= n\left[\frac{(1-\theta_0)\phi_0}{(1-\theta)^2} + \frac{\phi^2(1-\phi_0+\theta_0\phi_0)}{(1-\phi+\theta\phi)^2}\right] > 0.
\end{aligned}$$

This establishes that the KL divergence is strictly convex in $\theta$. Solving the critical equation $0 = dR/d\theta\,(\theta)$ or using the appropriate boundary point (for the monotonic cases) gives the minimising point $\theta = [(\phi - \phi_0 + \theta_0\phi_0)/\phi]_{0,1}$. Substituting $\theta = (\phi - \phi_0 + \theta_0\phi_0)/\phi$ gives:

$$\log\left(\frac{(1-\theta_0)\phi_0}{(1-\theta)\phi}\right) = \log\left(\frac{1-\phi_0+\theta_0\phi_0}{1-\phi+\theta\phi}\right) = 0,$$

which gives zero KL divergence. This establishes all the results to be shown. ■

**PROOF OF THEOREM 7:** This result is an application of the general posterior consistency results for Bayesian analysis set out in Schwarz (1965) and later expanded in LeCam (1973) and van der Vaart (1998). We will prove the result explicitly for the case where $\kappa = 0$ (i.e., for the left-HyperBeta distribution) so that $\theta_* = \theta_0\phi_0/\phi$. The corresponding proof for the case $\kappa = 1$ (i.e., for the right-HyperBeta distribution) is analogous.

The Schwartz theorem in Schwarz (1965) requires the model to meet two conditions — the "prior support" condition and the "testing condition". We now show that our model meets both of these conditions:

(a) **Prior support condition:** The present model uses the fixed parameter space $\Theta = [0,1]$ and a prior that encompasses the interior of the parameter space in its support. From Theorem 7A we know that the Kullback-Leibler (KL) divergence between $(\theta, \phi)$ and $(\theta_0, \phi_0)$ is:

$$\text{KL}(\theta_0, \phi_0 | \theta, \phi) = \sum_{s=0}^{n} \text{Bin}(s|n, \theta_0\phi_0) \begin{bmatrix} s \log\left(\dfrac{\theta_0\phi_0}{\theta\phi}\right) \\ +(n-s)\log\left(\dfrac{1-\theta_0\phi_0}{1-\theta\phi}\right) \end{bmatrix},$$

and this divergence is minimised uniquely at the point $\theta_* = \theta_*$ where the divergence is zero. Now, define the set $\Theta_0(\varepsilon) \equiv \{\theta | \text{KL}(\theta_0, \phi_0|\theta, \phi) \leq \varepsilon\}$ which is the parameter subspace for $\theta$ where the KL-divergence is no greater than $\varepsilon$ (i.e., this is the subspace of parameter values that are "close" to $\theta_*$ in the sense of having low KL divergence). Since the KL divergence is continuous in $\theta$ and $\text{KL}(\theta_0|\theta_*) = 0$ this means that for any $\varepsilon > 0$ there is a neighbourhood $\mathcal{U}_* \ni \theta_*$ with $\mathcal{U}_* \subseteq \Theta_0(\varepsilon)$. Since the prior distribution for the model encompasses the interior of the parameter space in its support we have:

$$\mathbb{P}(\theta \in \Theta_0(\varepsilon)) \geq \mathbb{P}(\theta \in \mathcal{U}_*) > 0,$$

which establishes the prior support condition for the theorem.

(b) **Testing condition:** The present model uses the fixed parameter space $\Theta = [0,1]$ which is compact and unidimensional. Moreover, the mapping $\theta \longmapsto f_\theta$ from the parameter to the sampling density is identifiable and continuous. Using Lemma 10.6 in van der Vaart (1998), it follows that the model satisfies the testing condition.

We have now established the prior support and testing conditions for the Schwartz theorem. The result for the left-HyperBeta distribution follows as a direct application of this theorem. The result for the right-HyperBeta distribution can be proved analogously. ■

**PROOF OF THEOREMS 8A AND 8B:** These results follow trivially from the definition of the functions $L_1$ and $R_1$ using the law of the unconscious statistician. ■

**PROOF OF THEOREM 9A:** We begin by showing the log-likelihood for a single observed value $x$ and then we extend this to the IID model in the theorem. The log-density of the left-HyperBeta distribution is given by:

$$\begin{aligned}\ell_x(\alpha,\beta,\gamma,\phi,0) &\equiv \log \text{LHyperBeta}(x|\alpha,\beta,\gamma,\phi) \\ &= (\alpha-1)\log(x) + (\beta-1)\log(1-x) + \gamma\log(1-x\phi) \\ &\quad - \log\left(\int_0^1 x^{\alpha-1}(1-x)^{\beta-1}(1-x\phi)^{\gamma}\,dx\right)\end{aligned}$$

This gives the partial derivatives:

$$\begin{aligned}\frac{\partial\ell_x}{\partial\alpha}(\alpha,\beta,\gamma,\phi,0) &= \log(x) - \frac{\partial}{\partial\alpha}\log\left(\int_0^1 x^{\alpha-1}(1-x)^{\beta-1}(1-x\phi)^{\gamma}\,dx\right) \\ &= \log(x) - \frac{\int_0^1 \log(x)\, x^{\alpha-1}(1-x)^{\beta-1}(1-x\phi)^{\gamma}\,dx}{\int_0^1 x^{\alpha-1}(1-x)^{\beta-1}(1-x\phi)^{\gamma}\,dx} \\ &= \log(x) - \int_0^1 \log(x)\,\text{LHyperBeta}(x|\alpha,\beta,\gamma,\phi)\,dx \\ &= \log(x) - L_1(\alpha,\beta,\gamma,\phi), \\ \frac{\partial\ell_x}{\partial\beta}(\alpha,\beta,\gamma,\phi,0) &= \log(1-x) - \frac{\partial}{\partial\beta}\log\left(\int_0^1 x^{\alpha-1}(1-x)^{\beta-1}(1-x\phi)^{\gamma}\,dx\right) \\ &= \log(1-x) - \frac{\int_0^1 \log(1-x)\, x^{\alpha-1}(1-x)^{\beta-1}(1-x\phi)^{\gamma}\,dx}{\int_0^1 x^{\alpha-1}(1-x)^{\beta-1}(1-x\phi)^{\gamma}\,dx} \\ &= \log(1-x) - \int_0^1 \log(1-x)\,\text{LHyperBeta}(x|\alpha,\beta,\gamma,\phi)\,dx \\ &= \log(1-x) - L_2(\alpha,\beta,\gamma,\phi), \\ \frac{\partial\ell_x}{\partial\gamma}(\alpha,\beta,\gamma,\phi,0) &= \log(1-x\phi) - \frac{\partial}{\partial\gamma}\log\left(\int_0^1 x^{\alpha-1}(1-x)^{\beta-1}(1-x\phi)^{\gamma}\,dx\right) \\ &= \log(1-x\phi) - \frac{\int_0^1 \log(1-x\phi)\, x^{\alpha-1}(1-x)^{\beta-1}(1-x\phi)^{\gamma}\,dx}{\int_0^1 x^{\alpha-1}(1-x)^{\beta-1}(1-x\phi)^{\gamma}\,dx} \\ &= \log(1-x\phi) - \int_0^1 \log(1-x\phi)\,\text{LHyperBeta}(x|\alpha,\beta,\gamma,\phi)\,dx \\ &= \log(1-x\phi) - L_3(\alpha,\beta,\gamma,\phi),\end{aligned}$$

$$\frac{\partial \ell_x}{\partial \phi}(\alpha,\beta,\gamma,\phi,0) = -\frac{\gamma x}{1-x\phi} - \frac{\partial}{\partial \phi}\log\left(\int_0^1 x^{\alpha-1}(1-x)^{\beta-1}(1-x\phi)^{\gamma}dx\right)$$

$$= -\frac{\gamma x}{1-x\phi} + \gamma \cdot \frac{\int_0^1 x^{\alpha}(1-x)^{\beta-1}(1-x\phi)^{\gamma-1}dx}{\int_0^1 x^{\alpha-1}(1-x)^{\beta-1}(1-x\phi)^{\gamma}dx}$$

$$= -\frac{\gamma x}{1-x\phi} + \gamma \cdot \frac{\text{B}(\alpha,\beta)\ \ {}_2F_1(-\gamma,\alpha,\alpha+\beta;\phi)}{\text{B}(\alpha+1,\beta)\ \ {}_2F_1(-\gamma-1,\alpha+1,\alpha+\beta+1;\phi)}.$$

The resulting log-likelihood and partial derivatives for the score function follow by summing over the observations $x_1, \dots, x_n$ used in the theorem, which establishes the result. ■

**PROOF OF THEOREM 9B:** We begin by showing the log-likelihood for a single observed value $x$ and then we extend this to the IID model in the theorem. The log-density of the right-HyperBeta distribution is given by:

$$\begin{aligned}\ell_x(\alpha,\beta,\gamma,\phi,1) &\equiv \log \text{RHyperBeta}(x|\alpha,\beta,\gamma,\phi) \\ &= (\alpha-1)\log(x) + (\beta-1)\log(1-x) + \gamma\log(1-\phi+x\phi) \\ &\quad -\log\left(\int_0^1 x^{\alpha-1}(1-x)^{\beta-1}(1-\phi+x\phi)^{\gamma}dx\right)\end{aligned}$$

This gives the partial derivatives:

$$\begin{aligned}\frac{\partial \ell_x}{\partial \alpha}(\alpha,\beta,\gamma,\phi,1) &= \log(x) - \frac{\partial}{\partial \alpha}\log\left(\int_0^1 x^{\alpha-1}(1-x)^{\beta-1}(1-\phi+x\phi)^{\gamma}dx\right) \\ &= \log(x) - \frac{\int_0^1 \log(x)\, x^{\alpha-1}(1-x)^{\beta-1}(1-\phi+x\phi)^{\gamma}dx}{\int_0^1 x^{\alpha-1}(1-x)^{\beta-1}(1-\phi+x\phi)^{\gamma}dx} \\ &= \log(x) - \int_0^1 \log(x)\,\text{RHyperBeta}(x|\alpha,\beta,\gamma,\phi)\,dx \\ &= \log(x) - R_1(\alpha,\beta,\gamma,\phi),\end{aligned}$$

$$\begin{aligned}\frac{\partial \ell_x}{\partial \beta}(\alpha,\beta,\gamma,\phi,1) &= \log(1-x) - \frac{\partial}{\partial \beta}\log\left(\int_0^1 x^{\alpha-1}(1-x)^{\beta-1}(1-\phi+x\phi)^{\gamma}dx\right) \\ &= \log(1-x) - \frac{\int_0^1 \log(1-x)\, x^{\alpha-1}(1-x)^{\beta-1}(1-\phi+x\phi)^{\gamma}dx}{\int_0^1 x^{\alpha-1}(1-x)^{\beta-1}(1-\phi+x\phi)^{\gamma}dx}\end{aligned}$$

$$= \log(1-x) - \int_0^1 \log(1-x)\,\mathrm{RHyperBeta}(x|\alpha,\beta,\gamma,\phi)\,dx$$

$$= \log(1-x) - R_2(\alpha,\beta,\gamma,\phi),$$

$$\frac{\partial \ell_x}{\partial \gamma}(\alpha,\beta,\gamma,\phi,1) = \log(1-\phi+x\phi) - \frac{\partial}{\partial \gamma}\log\left(\int_0^1 x^{\alpha-1}(1-x)^{\beta-1}(1-\phi+x\phi)^{\gamma}dx\right)$$

$$= \log(1-\phi+x\phi) - \frac{\int_0^1 \log(1-\phi+x\phi)\,x^{\alpha-1}(1-x)^{\beta-1}(1-\phi+x\phi)^{\gamma}dx}{\int_0^1 x^{\alpha-1}(1-x)^{\beta-1}(1-\phi+x\phi)^{\gamma}dx}$$

$$= \log(1-\phi+x\phi) - \int_0^1 \log(1-\phi+x\phi)\,\mathrm{RHyperBeta}(x|\alpha,\beta,\gamma,\phi)\,dx$$

$$= \log(1-\phi+x\phi) - R_3(\alpha,\beta,\gamma,\phi),$$

$$\frac{\partial \ell_x}{\partial \phi}(\alpha,\beta,\gamma,\phi,1) = -\frac{\gamma(1-x)}{1-\phi+x\phi} - \frac{\partial}{\partial \phi}\log\left(\int_0^1 x^{\alpha-1}(1-x)^{\beta-1}(1-\phi+x\phi)^{\gamma}dx\right)$$

$$= -\frac{\gamma(1-x)}{1-\phi+x\phi} + \gamma\cdot\frac{\int_0^1 x^{\alpha-1}(1-x)^{\beta}(1-\phi+x\phi)^{\gamma-1}dx}{\int_0^1 x^{\alpha-1}(1-x)^{\beta-1}(1-\phi+x\phi)^{\gamma}dx}$$

$$= -\frac{\gamma(1-x)}{1-\phi+x\phi} + \gamma\cdot\frac{\mathrm{B}(\alpha,\beta)\ {}_2F_1(-\gamma,\alpha,\alpha+\beta;\phi)}{\mathrm{B}(\alpha,\beta+1)\ {}_2F_1(-\gamma-1,\alpha,\alpha+\beta+1;\phi)}.$$

The resulting log-likelihood and partial derivatives for the score function follow by summing over the observations $x_1, \dots, x_n$ used in the theorem, which establishes the result. ■